\documentclass[12pt,a4paper]{article}
\usepackage[T1]{fontenc}
\usepackage{lmodern}
\usepackage[latin1]{inputenc}
\usepackage{latexsym}                 
\usepackage{color}                 
\usepackage{amssymb}
\usepackage{amsmath}
\usepackage{amsthm}
\usepackage{graphicx}
\usepackage{hyperref}
\usepackage{mathptmx}      % use Times fonts if available on your TeX system
\def\RSlabel#1{\label{#1}%
}
\newcommand{\bql}[1]{%
\begin{equation}\label{#1}%
}
\def\RScite#1{\cite{#1}%
}
\def\RSref#1{\ref{#1}%
}
\def\eref#1{(\ref{#1}%
)}
\newtheorem{Theorem}{Theorem}
\newtheorem{Definition}{Definition}
\def\calh{{\cal H}}
\def\fa{\hbox{ for all }}
\def\biglf{\par\bigskip\noindent}
\newcommand{\ep}[1]{\ensuremath{\varepsilon^{#1}}}
\newcommand{\R}{\ensuremath{\mathbb{R}}}
\newcommand{\eq}{\end{equation}}
\begin{document}
%%%%%%%%%%%%%%%%%%%%%%%%%%%%%%%%%%%%%%%%%%%%%%%%%%%%%%%%%
\title{Scaling of Radial Basis Functions}
%\renewcommand*\thefootnote{\alph{footnote}}
%\author{E. Larsson\footnote{Univ. of Uppsala, Sweden}
%  and R. Schaback\footnote{Institut f\"ur Numerische und Angewandte
%Mathematik, Universit\"at G\"ottingen, Lotzestra\ss{}e 16--18,
%D--37083 G\"ottingen, Germany}}
%\footnotetext[1]{Univ. of Uppsala, Sweden}
%\footnotetext[2]{Institut f\"ur Numerische und Angewandte
%Mathematik, Universit\"at G\"ottingen, Lotzestra\ss{}e 16--18,
%D--37083 G\"ottingen, Germany}
\author{E. Larsson\phantom{\footnotemark[1]}$^{1,}$
  \phantom{\footnotemark[3]}$^3$
  and R. Schaback\phantom{\footnotemark[2]}$^{2,}$
  \phantom{\footnotemark[3]}$^3$
}
\footnotetext[1]{Scientific Computing, Department of Information Technology, Uppsala University, Box 337, SE-751 05 Uppsala, Sweden}
\footnotetext[2]{Institut f\"ur Numerische und Angewandte
Mathematik, Universit\"at G\"ottingen, Lotzestra\ss{}e 16--18,
D--37083 G\"ottingen, Germany}
\footnotetext[3]{Authors in alphabetical order}
\maketitle
\abstract{This paper studies the influence of
  scaling on the behavior of Radial Basis Function interpolation.
  It focuses on
  certain central aspects, but does not try to be exhaustive.
  The most important questions are: How does the error of a
  kernel-based interpolant
  vary with the scale of the kernel chosen? How does the standard
  error bound vary? And since fixed functions may be in spaces that
  allow scalings, like global Sobolev spaces, is there a scale of
  the space that matches the function best?
  The last question is answered in the affirmative for Sobolev spaces,
  but the required scale may be hard to estimate. Scalability of functions
  turns out to be restricted for spaces generated by analytic kernels,
  unless the functions are band-limited. In contrast to other papers,
  polynomials and polyharmonics
  are included as {\em flat limits}  when checking scales experimentally,
  with an independent computation. The numerical results show that
  the hunt for near-flat scales is questionable, if users
  include the flat limit cases right from the start.
  When  there are not enough data to evaluate errors directly,
  the scale of the standard error bound can be varied, up to replacing the norm
  of the unknown function by the norm of the interpolant.
  This follows the behavior of the actual error qualitatively well,
  but is only of limited value for estimating error-optimal scales.
  For kernels and functions with unlimited smoothness, the given interpolation
  data are proven to be insufficient for determining useful scales. 
  
  %***************************
%****************************************************************
\section{Introduction}\RSlabel{SecIntro}
Throughout the paper it is assumed that readers are familiar
with the basics of kernel-based methods, e.g. from books by M.D. Buhmann
\cite{buhmann:2003-1}, H. Wendland \cite{wendland:2005-1},
and G. Fasshauer/M. McCourt \cite{fasshauer-mccourt:2015-1}.
\biglf
It is a well-known fact 
% (see e.g. \cite{Schaback-wendland:2006-1}) 
that interpolation and approximation with radial
basis functions is
crucially dependent on {\em scaling}. There are two ways to
introduce scaling into a radial basis function. For instance, the classical
multiquadrics are often written as kernels
$$
\Phi(\|x-y\|_2):=(c^2+\|x-y\|_2^2)^\beta=c^{2\beta}(1+\|x-y\|_2^2/c^2)^\beta
\fa x,y\in\R^d,
$$
and compactly supported Wendland functions are scaled like in 
$$
\Phi(\|x-y\|_2):=(1-\|x-y\|_2/c)_+^m\;p_{d,k}(\|x-y\|_2/c)\fa x,y\in\R^d
$$
to let the support radius be $c>0$. The PDE--oriented literature
often calls $c$ a {\em shape parameter}
(see e.g. \cite{kansa-carlson:1992-1}). In 
another context started by B. Fornberg and his collaborators
\cite{driscoll-fornberg:2002-1} 
one works with $\epsilon=1/c$ instead
and calls $\epsilon\to 0$ the {\em flat limit}.
For details on flat limits, see
e.g. \RScite{driscoll-fornberg:2002-1,fornberg-et-al:2004-1,%
  larsson-fornberg:2005-1,schaback:2005-2,schaback:2008-2,song-et-al:2012}
with a good summary in \RScite{fasshauer-mccourt:2015-1}.
Here, we shall 
stick to the latter notation and and consider scaled radial basis functions    
$$
\Phi_\epsilon(\|x-y\|_2):=\Phi(\epsilon\|x-y\|_2)
$$ 
for $x,\,y\in\R^d$ with $\Phi\;:\;[0,\infty)\to\R$.
We ignore multiplication of the kernel by a scalar here, because the scalar
cancels out for interpolation and approximation. However, it turns up
as a {\em process variance}
in methods that apply probabilistic estimation after rewriting
kernel methods via Gaussian processes. See
\cite{fasshauer-mccourt:2015-1,scheuerer-et-al:2013-1}
for the connection to the deterministic situation.
\biglf
If data $f(x_1),\ldots,f(x_n)\in\R$ of a function $f$
are to be interpolated in scattered locations
$x_1,\ldots,x_n\in\R^d$, they determine coefficients $\alpha_j$ by the
linear system
\bql{eqsys}
f(x_k)=\displaystyle{\sum_{j=1}^n\alpha_j \Phi_\epsilon(x_k,x_j) ,\;
1\leq k\leq n  }
\eq
with a positive definite {\em kernel matrix} $A^{X,\Phi_\epsilon}$
with entries 
$\Phi_\epsilon(x_j,x_k),\;1\leq j,k\leq n$.
Then
the interpolant is written as
$$
s_{f,X,\Phi_\epsilon}(x)=\displaystyle{\sum_{j=1}^n\alpha_j
  \Phi_\epsilon(x,x_j)   },
$$
or as a linear combination
$$
s_{f,X,\Phi_\epsilon}(x)=\displaystyle{\sum_{j=1}^nf(x_j)u_{x_j}^{X,\Phi_\epsilon}(x)}
$$
for {\em Lagrangians} or {\em cardinal interpolants}
with $u_{x_j}^{X,\Phi_\epsilon}(x_{i})=\delta_{ji}$
given by 
$$
u_{x_j}^{X,\Phi_\epsilon}(x)=\displaystyle{
  \sum_{k=1}^n \beta^{X,\Phi_\epsilon}_{jk}\Phi_\epsilon(x,x_k)  } 
$$
via the coefficients $\beta_{jk}^{X,\Phi_\epsilon},\;1\leq j,k\leq n$ from
the inverse
of the kernel matrix $A^{X,\Phi_\epsilon}$ because of
  $$
u_{x_j}^{X,\Phi_\epsilon}(x_i)=\displaystyle{\delta_{ji}=
  \sum_{k}=1}^n \beta^{X,\Phi_\epsilon}_{jk}\Phi_\epsilon(x_i,x_k),
1\leq j,i\leq n.
  $$
Theoretically, one can use arbitrary positive definite
radial basis functions with arbitrary scales. The result can be seen as
independent of $f$, if one only looks at the Lagrangians.
This will be the topic of Section \RSref{SecSoL}.
\biglf
Fixing a norm for $f$ and the interpolant,
one can ask for an {\em error-optimal} scale that minimizes
$\|f-s_{f,X,\Phi_\epsilon}\|$ over all manageable scales. 
This is an interesting
problem in theory and practice, and still open in both respects. There
is a vast literature on numerical methods % \red{citations}
for estimating optimal scales in practice,
but we shall not compare these here. 
\biglf
If users minimize errors or error bounds for a specific $X$, the
distribution and density of points may have a strong influence on the result,
leading to different suggested
scales for different $X$. Users should be cautious when seeing
estimates of optimal scales behaving like $1/h(X,\Omega)$
for the 
 {\em fill distance}
$$
h=h(X,\Omega)=\displaystyle{\sup_{y\in \Omega}\min_{x\in X}\|x-y\|_2,   } 
$$
since it is well-known \cite{buhmann:1989-2,baxter-brummelhuis:2022-1}
that there is no convergence in such a
{\em stationary} situation except for
  conditionally positive kernels. 
\biglf
A possible workaround to minimizing the error directly
is to ask for an {\em error-bound-optimal} scale
that minimizes some specified norm of the error.
The standard pointwise $L_\infty$ error bound is
\bql{eqfsPfpw}
|f(x)-s_{f,X,\Phi_\epsilon}(x)|\leq P_{X,\Phi_\epsilon}(x)
\|f\|_{\Omega,\Phi_\epsilon}
\eq
using the domain-independent {\em Power Function}. 
It holds for all functions $f$ from the local {\em native space }
  $\calh_{\Omega,\Phi}$  and all finite sets $X\subset\Omega$ of centers.
  By standard extension arguments \cite[Section 10.7]{wendland:2005-1} one has
  $$
\|f\|_{\Omega,\Phi_\epsilon}\leq \|E_\Omega(f)\|_{\R^d,\Phi_\epsilon}
$$
for the canonical extension $E_\Omega(f)$ of $f$ to the whole space.
If norms are taken, one gets two variations. The local bound is
\bql{eqfsPfloc}
\|f-s_{f,X,\Phi_\epsilon}\|_{\Omega,\infty}\leq \|P_{X,\Phi_\epsilon}\|_{\Omega,\infty}
\|f\|_{\Omega,\Phi_\epsilon}\leq \|P_{X,\Phi_\epsilon}\|_{\Omega,\infty}
\|E_\Omega f\|_{\R^d,\Phi_\epsilon}
\eq
holding on the local native space $\calh_{\Omega,\Phi}$,
while a global one is
\bql{eqfsPf}
\|f-s_{f,X,\Phi_\epsilon}\|_{\Omega,\infty}\leq \|P_{X,\Phi_\epsilon}\|_{\Omega,\infty}
\|f\|_{\R^d,\Phi_\epsilon}
\eq
holding on the global native space $\calh_{\R^d,\Phi}$.
This is the version we
shall analyze in Section \RSref{SecSoEB}.
Note that \eref{eqfsPf} turns into
\eref{eqfsPfloc} if $f$ is replaced by $E_\Omega(f)$. 
Both factors in the right-hand side of \eref{eqfsPf} will depend on
$\epsilon$. The bound splits into an $f$-dependent part
$\|f\|_{\R^d,\Phi_\epsilon}$ and an $X$-dependent part
$\|P_{X,\Phi_\epsilon}\|_{\Omega,\infty}$ whatever the scaling is. This raises the
question whether there is a scale $\epsilon$ that minimizes
$\|f\|_{\R^d,\Phi_\epsilon}$ without caring for the data.
\biglf
Furthermore,
convergence rates for interpolation
often do not depend on scaling, but the factor in front of the rates
does. This means that there must be an $X$-independent part in the whole
problem.
\biglf
All of this raises the question whether 
functions have a ``natural scaling'' that can be recovered approximately
somehow, being independent of the interpolation problem. And since
$\|f\|_{\Phi_\epsilon}$ is a factor in the error bound \eref{eqfsPf},
there is some hope that a minimal $\|f\|_{\Phi_\epsilon}$
also cares for small factors in front of convergence rates.
Surprisingly, functions in global Sobolev spaces have such a
natural scale,
as proven in Section \RSref{SecNSoF}. 
\biglf
But there are practical limits to scaling.
It is well--known \cite{schaback:1995-1}
that the condition of the kernel matrices $A^{X,\Phi_\epsilon}$
increases dramatically
for $\epsilon \to 0$, while the interpolants still exist for all $\epsilon>0$.
There are many workarounds for this, e.g. Contour-Pad\'e 
\cite{fornberg-wright:2004-1}, RBF-QR \cite{fornberg-et-al:2011-1}, and 
RBF-GA \cite{fornberg-et-al:2013-1} by the group around Bengt Fornberg,
and Hilbert-Schmidt-SVD by Fasshauer/McCourt
\cite[Chapter 13]{fasshauer-mccourt:2015-1}. 
But we shall
only focus on the flat limit here, not on methods to reach it via small scales.
It requires to 
distinguish between two kinds  of radial kernels:
the {\em analytic} ones have convergent
expansions into powers of $\|x-y\|_2^2$,
like the Gaussian or inverse multiquadrics $\Phi(x-y)=(1+\|x-y\|_2^2)^{-m}$
for positive $m$. These have infinite smoothness, and their Fourier transform
decays exponentially towards infinity. The other
{\em non-analytic} kind has limited smoothness,
and their Fourier transform
decays algebraically, i.e. as  a negative finite power towards infinity.
These include Mat\'ern-Sobolev and compactly supported Wendland
\RScite{wendland:1995-1}
or Buhmann \RScite{buhmann:1998-1} kernels,
among others.
The flat limit for the analytic kernels is a polynomial
\cite{driscoll-fornberg:2002-1} 
except
for certain degenerations \cite{larsson-fornberg:2005-1,
  lee-et-al:2005-1,schaback:2008-2}
depending on the point sets, while the flat limit
for the others \RScite{song-et-al:2012}
is a polyharmonic interpolant based on the conditionally positive
definite kernels 
$$
K(x,y)=\left\{
\begin{array}{ll}
\|x-y\|^{2m-d} & d \hbox{ odd}\\ 
\|x-y\|^{2m}\log \|x-y\|_2 & d \hbox{ even}\\ 
\end{array}
\right.
$$
with generalized Fourier transforms $\|\omega\|_2^{-2m}$ on $\R^d$.
These kernels
are scale-invariant, because the scale comes to the front as a scalar
factor in the Fourier transform.
They are hard to beat if it comes to
handle stencils for approximation of derivatives
\RScite{davydov-schaback:2019-1}, and they offer a convenient
bypass around all scaling problems. 
Situations with very small $\epsilon$ and
kernels with finite smoothness will be explained perfectly
by interpolation with polyharmonic kernels
that should be used right from the
start in such cases. The discussion of flat limits should
therefore be confined
to analytic kernels. 
\biglf
Summarizing, the above discussion shows that
the strong dependence of radial basis 
function techniques on scaling is a feature, not a bug.
The functions supplying the data already have a hidden natural scale,
independent of how the reconstruction by interpolation
or approximation is done, and good recovery methods
should therefore not be scale--independent.
This paper tries to clarify the scaling effects to some extent.
\biglf
It starts by collecting some basic facts on scaling in Section
\RSref{SecBFSaS} for the convenience of readers. These include
some useful invariance
relations and describe the dependence of Power Functions and
Lagrangians on scaling. Optimal scales of functions result from
studying norms $\|f\|_{\Phi_\epsilon}$ for fixed $f$ as functions
of $\epsilon$ in section \RSref{SecSoKbN}.
\biglf
Finally, we look at the limit $\epsilon\to 0$ for analytic kernels $\Phi$
in section \RSref{SecEE}, ignoring possible degenerations.
It is known that the limit interpolant is a polynomial, 
but here we study the behavior of the norm of the
interpolation error
as a function of $\epsilon$. Experimentally, it is a
function of $\epsilon$  that can have sharp local minima at
seemingly unexpected scales,
but there are also many cases where the flat limit
has an optimal error norm.
Certain criteria for these two cases are provided,
but they are hard to handle in practice.
However, it is proven that
any fixed set of interpolation data does not
determine whether the flat limit is
optimal or not. Users need additional data for error evaluation
when they search for optimal scales.
%****************************************************************
\section{Basic Facts on Scaling and Spaces}\RSlabel{SecBFSaS}
Throughout the paper, functions will be real-valued, defined on $\R^d$,
and Fourier transformable. We do not treat localized versions of scaling here.
Furthermore, we focus on functions with continuous point evaluation,
and therefore we work on Hilbert spaces $\calh_\Phi$
with reproducing kernels $\Phi$.
These kernels will be translation-invariant and Fourier transformable, i.e.
there is a reproduction property
$$
f(x)=(f,\Phi(x-\cdot))_{\calh_\Phi} \fa f\in \calh_\Phi,\;x\in\R^d
$$
and an inner product
$$
(f,g)_{\calh_\Phi}=(2\pi)^{-d/2}\int_{\R^d}\dfrac{\hat f(\omega)\overline{\hat g(\omega)}}{\hat
  \Phi(\omega)}d\omega \fa f\in{\calh_\Phi}
$$
with the usual property
\bql{eqKKrep}
\Phi(x-y)=(\Phi(x-\cdot),\Phi(\cdot-y))_{\calh_\Phi}\fa x,y\in\R^d.
\eq
A particularly interesting case is
Sobolev space $W_2^m(\R^d)$ for $m>d/2$ with the
Whittle-Mat\'ern kernel
\bql{eqWMK}
\Phi(x):=\dfrac{2^{1-m}}{\Gamma(m)}\|x\|_2^{m-d/2}K_{m-d/2}(\|x\|_2)  
\eq
using the modified Bessel function of second kind.
It has the $d$-variate Fourier transform
\bql{eqWMKFT}
\hat \Phi(\omega)=(1+\|\omega\|_2^2)^{-m}\fa \omega\in\R^d
\eq
to make it compatible with the inner product above.
\biglf
\begin{Definition}\RSlabel{DefScaling}
  A function $f\;:\;\R^d\to \R$ will be {\em scaled} by
 \bql{eqscale}
%RS f_\epsilon(x):=\epsilon^{d/2}\;f(\epsilon x) \fa x\in\R^d,\,\epsilon>0.
f_\epsilon(x):=f(\epsilon x) \fa x\in\R^d,\,\epsilon>0.
\eq
\end{Definition}
\noindent
Note that scaling is ``inverted'' in the Fourier domain by
\bql{eqFTI}
\widehat{f_\epsilon} (\omega)
=
\epsilon^{-d}\hat f (\omega/\epsilon)
=
\epsilon^{-d}\hat f_{1/\epsilon} (\omega)\fa \omega\in\R^d,\;\epsilon>0.
\eq
Thus it is equivalent up to a factor to consider scaling of a function or its
Fourier transform. If one would define scaling differently, namely by
$f_\epsilon(x):=\epsilon^{d/2}f(\epsilon x)$, there would be Fourier transform
symmetry. We avoid this, because we want to keep the fact that $\epsilon\to 0$
implies that $f_\epsilon$ does not go to zero,
while the Fourier transform of $f_\epsilon$ goes to zero if it decays fast
enough at infinity. This is the ``flat limit'' situation, and 
this paper will keep an eye on flat limits throughout.
\biglf
The scaling law \eref{eqscale} will scale compact supports properly,
keeping $f_\epsilon(0)=f(0)$ invariant, but this will not be true
in frequency space. There, the integral over frequency space is invariant.
%****************************************************************
% \subsection{Scaling of Kernels }\RSlabel{SecSoK}
\biglf
When we scale the kernel $K$ of the Hilbert space ${\calh_\Phi}$, we
shall denote the scaled inner product by $(\cdot,\cdot)_{\Phi_\epsilon}$
belonging to the scaled
kernel $K_\epsilon$ of the Hilbert space ${\calh_{\Phi_\epsilon}}$
Then the dual version
$$
(\delta_x,\delta_y)_\Phi=(\Phi(x-\cdot),\Phi(\cdot-y))_{\calh_\Phi}
=\Phi(x-y) \fa x,y\in\R^d
$$
of \eref{eqKKrep} yields
\bql{eqscaleddeltas}
(\delta_x,\delta_y)_{\Phi_\epsilon}=\Phi_\epsilon(x-y)
=\Phi(\epsilon x -\epsilon y)=
(\delta_{\epsilon x},\delta_{\epsilon  y})_{\Phi}
\eq
for all $ x,y\in \R^d,\;\epsilon >0$. This is a {\em scaling law}
for point evaluation functionals.
%****************************************************************
%\subsection{Nesting of Kernel-Based Spaces}\RSlabel{SecNoKbS}
\biglf
When treating fixed functions $f\in \calh_\Phi$
with scaled kernels $K_{\Phi_\epsilon}$,
it is not clear whether all scaled native spaces  $\calh_{\Phi_\epsilon}$
contain $f$, and whether the native spaces are nested or norm-equivalent.
We postpone this to the study of norms $\|f\|_{\Phi_\epsilon}$ for
$f\in \calh_{\Phi}$, and it will turn out that analytic kernels
will cause problems.
%****************************************************************
\subsection{Scaling of Kernel-Based Norms }\RSlabel{SecSoKbN}
We consider functions $f\in \calh_\Phi$ and start with a simple scaling law
\bql{eqfscaleinv}
\|f_\epsilon\|_{\Phi_\epsilon}=\|f\|_{\Phi}
\eq
that follows from
$$
\begin{array}{rcl}
  \|f_\epsilon\|_{\Phi_\epsilon}^2
  &=&
\displaystyle{ \int_{\R^d}
  \dfrac{|\hat f_\epsilon(\omega)|^2}{\hat \Phi_\epsilon(\omega)}d\omega}\\
  &=&
  \epsilon^{-2d}\displaystyle{ \int_{\R^d} \dfrac{|\hat f(\omega/\epsilon)|^2}{\hat \Phi_\epsilon(\omega)}d\omega}\\
  &=&
  \epsilon^{-d}\displaystyle{ \int_{\R^d} \dfrac{|\hat f(\eta)|^2}%
    {\hat \Phi_\epsilon(\eta\epsilon)}d\eta}\\
  &=&
  \displaystyle{ \int_{\R^d} \dfrac{|\hat f(\eta)|^2}%
    {\hat \Phi(\eta)}d\eta}\\
  &=&
  \|f\|^2_{\Phi}.
\end{array}
$$
In the form
$$
\|f_\epsilon\|_\Phi = \|f\|_{\Phi_{1/\epsilon}}
$$
this proves that scaling a function or a kernel is the same thing
for calculating native space norms,
as long as one of the sides exists.
\biglf
Now we fix functions $f\in\calh_\Phi$ and check for which kernel scales
we have $f\in\calh_{\Phi_\epsilon}$. 
A scale $\epsilon$ is called $\Phi$-{\em admissible for} $f\in\calh_\Phi$,
if $\|f\|_{\Phi_\epsilon}$ is finite.
\biglf
A simple restriction on admissible scales is
\bql{eqfPlowbnd}
\begin{array}{rcl}
  \|f\|_{\Phi_\epsilon}^2
  &=& \epsilon^{d}\int_{\R^d}
  \dfrac{|\hat f(\omega)|^2}{\hat \Phi(\omega/\epsilon)}d\omega\\
  &\geq &
  \epsilon^{d}\|f\|_{L_2}\|\hat \Phi\|_\infty^{-1}
\end{array}
\eq
with three implications: for large $\epsilon$ the norm  $\|f\|_{\Phi_\epsilon}^2$
must be large, for small  $\epsilon$ it can not be smaller than
${\cal O}(\epsilon^d)$, and admissible scales $\epsilon$ with
$\|f\|_{\Phi_\epsilon}^2\leq \|f\|_{\Phi}^2$ are bounded above.
Altogether, the case of large $\epsilon$ is not interesting. 
\begin{Theorem}\RSlabel{TheadmKerScaFin}
For kernels with finite smoothness, arbitrary scales are admissible.
In particular, if the Fourier transform of $\Phi$ behaves like
$\hat\Phi(\omega) =\Theta(\|\omega\|^{-\beta})$ with $\beta>d$
near infinity, then
$$
\|f\|_{\Phi_\epsilon}^2=\Theta(\epsilon^d 
\max(1,\epsilon^{-\beta}))\|f\|_{\Phi}^2
$$
for all $\epsilon>0$ and all $f\in \calh_\Phi$. The spaces
$\calh_{\Phi_\epsilon}$ are identical as sets and norm-equivalent,
the equivalence constants behaving like $\epsilon^d \max(1,\epsilon^{-\beta})$
\biglf
For Sobolev space
$W_2^ m(\R^d)$ with integer $m$, this holds for $\beta=2m$, but
we can get the explicit formula
 $$
\|f\|^2_{\Phi_\epsilon}=\epsilon^d\displaystyle{ 
\sum_{j=0}^m{m\choose j}\epsilon^{-2j}
|f|^2_{W_2^j(\R^d)}}\fa f\in W_2^m(\R^d).
 $$
Since we are in global Sobolev space, the seminorms are norms,
and thus none of the $|f|_{W_2^j(\R^d)}$  can vanish.
\end{Theorem}
\begin{proof}
We consider
\bql{eqfPeiRd}
\begin{array}{rcl}
  \|f\|_{\Phi_\epsilon}^2
  &=& \int_{\R^d}
  \dfrac{|\hat f(\omega)|^2}{\hat \Phi_\epsilon(\omega)}d\omega\\
  &=& \epsilon^{d}\int_{\R^d}
  \dfrac{|\hat f(\omega)|^2}{\hat \Phi(\omega/\epsilon)}\\
  &=& \epsilon^{d}\int_{\R^d}
  \dfrac{|\hat f(\omega)|^2}{\hat \Phi(\omega)}
        \dfrac{\hat \Phi(\omega)}{\hat \Phi(\omega/\epsilon)}
        d\omega\\
\end{array}
\eq
and get the two-sided bound
$$ %\bql{eqtwoside}
\inf_{\omega\in\R^d}\dfrac{\hat \Phi(\omega)}{\hat \Phi(\omega/\epsilon)}  
\leq
\dfrac{\|f\|_{\Phi_\epsilon}^2}{\epsilon^{d} \|f\|_{\Phi}^2 }
\leq 
\sup_{\omega\in\R^d}
        \dfrac{\hat \Phi(\omega)}{\hat \Phi(\omega/\epsilon)}.  
$$ %\eq
To show that both bounds behave like $\max(1,\epsilon^{-\beta})$
for the finite smoothness case,     
we assume
$$
\begin{array}{rclrclcll}
0&<&C_0&\leq & \hat\Phi(\omega) &\leq& C_1 &\hbox{for }\|\omega\|_2\leq c\\   
0&<&c_0\|\omega\|_2^{-\beta}&\leq & \hat\Phi(\omega)
&\leq& c_1\|\omega\|_2^{-\beta} &\hbox{for }\|\omega\|_2\geq c\\   
\end{array} 
$$
and bound the quotient by 
$$
\begin{array}{rclrcll}
\dfrac{C_0}{C_1}&\leq &\dfrac{\hat \Phi(\omega)}{\hat \Phi(\omega/\epsilon)}
&\leq &\dfrac{C_1}{C_0} &
\hbox{for } \|\omega\|_2\leq c, \|\omega\|_2/\epsilon \leq c,\\  
\dfrac{C_0c^\beta}{c_1}&\leq &\dfrac{\hat \Phi(\omega)}{\hat \Phi(\omega/\epsilon)}
&\leq &\epsilon^{-\beta}\dfrac{C_1c^\beta}{c_0} &
\hbox{for } \|\omega\|_2\leq c, \|\omega\|_2/\epsilon \geq c,\\  
\epsilon^{-\beta}\dfrac{c_0c^\beta}{C_1}&\leq &\dfrac{\hat \Phi(\omega)}{\hat \Phi(\omega/\epsilon)}
&\leq &\dfrac{c_1c^{-\beta}}{C_0} &
\hbox{for } \|\omega\|_2\geq c, \|\omega\|_2/\epsilon \leq c,\\  
\epsilon^{-\beta} \dfrac{c_0}{c_1}&\leq &\dfrac{\hat \Phi(\omega)}{\hat \Phi(\omega/\epsilon)}
&\leq &\epsilon^{-\beta} \dfrac{c_1}{c_0}&
\hbox{for } \|\omega\|_2\geq c, \|\omega\|_2/\epsilon\geq c.\\  
\end{array}
$$
If $m$ is an integer, we can continue from \eref{eqfPeiRd} to
\bql{eqSobEpsExact}
\begin{array}{rcl}
\|f\|^2_{\Phi_\epsilon}
&=&
\displaystyle{ (2\pi)^{-d/2}
\epsilon^d\sum_{j=0}^m{m\choose j}\epsilon^{-2j}\int_{\R^d}
|\hat{f}(\omega)|^2 \|\omega\|_2^{2j}
d\omega}\\
&=&
\displaystyle{ 
\epsilon^d\sum_{j=0}^m{m\choose j}\epsilon^{-2j}
|f|^2_{W_2^j(\R^d)}}.
\end{array} 
\eq
\end{proof}
\begin{Theorem}\RSlabel{TheadmKerScaInf}
For Gaussians and other kernels with exponential decay of the Fourier transform
at infinity, the inclusion $\calh_\Phi\subseteq\calh_{\Phi_\epsilon}$
is true only for $\epsilon\geq 1$. For single functions,
kernel scales $\epsilon<1$ may be admissible, but this depends on the function.
For bandlimited functions and all kernels, all kernel scales are admissible. 
\end{Theorem} 
\begin{proof}
If we assume an exponential law
$$
\hat\Phi(\omega)=c\exp(-\gamma\|\omega\|_2)
$$
near infinity, we get
$$
\begin{array}{rcl}
  \dfrac{\hat \Phi(\omega)}{\hat \Phi(\omega/\epsilon)}
  &=&
  \exp(-\gamma\|\omega\|_2)\exp(\gamma\|\omega\|_2/\epsilon).  
\end{array}
$$
Now $\|f\|^2_{\Phi_\epsilon}$ gets unbounded for $\epsilon<1$
provided that the integrable function
$|\hat f(\omega)|^2/\Phi(\omega)$ is bounded below near infinity 
by an arbitrarily large negative power of $\|\omega\|_2$.
The same argument works for the Gaussian, and shows that
admissible scales are strongly
$f$-dependent.
\biglf
For bandlimited functions, all scales are admissible by the above argumentation.
If the spectrum of $f$ is limited by $\|\omega\|_2\leq B$,
and if we define
  $$
\displaystyle{ \delta^-_\Phi(K):=\inf_{\|\omega\|_2\leq K}\hat\Phi(\omega)},\;
\displaystyle{ \delta^+_\Phi(K):=\sup_{\|\omega\|_2\leq K}\hat\Phi(\omega)},
  $$
equation \eref{eqfPeiRd} yields
\bql{eqfeffe}
\begin{array}{rclcl}
  \displaystyle{\dfrac{\delta^-_\Phi(B)}{\delta^+_\Phi(B/\epsilon)} } 
&\leq&
\dfrac{\|f\|_{\Phi_\epsilon}^2}{\epsilon^{d} \|f\|_{\Phi}^2 }
&\leq& 
\displaystyle{
        \dfrac{\delta^+_\Phi(B)}{\delta^-_\Phi(B/\epsilon)}}\\  
\end{array} 
\eq
to show admissibility of all scales for all kernels.
\end{proof}
Because the native space norm does not exist in the limit,
Theorem \RSref{TheadmKerScaInf} forbids to go 
to flat limits for analytic kernels
if $\|f\|^2_{\Phi_\epsilon}$ occurs, e.g. for standard
error bounds in terms of the Power Function. The flat limit of Lagrangians
in the left-hand side of \eref{eqphiscaleonly}
  exists, but the right-hand side cannot be formed for $\epsilon<1$.
\biglf
Here, we ignored domain-dependent native spaces and their norms.
  These are always equivalent to the global norms on $\R^d$,
  but the norm equivalence constants will vary when kernels are
  scaled. Since Section \RSref{SecSoEB} will be confined to
  global norms as well, we ignore the localized norms.
%****************************************************************
\subsection{Scaling the Power Function}\RSlabel{SecStPF}
For interpolation in $X=\{x_1,\ldots,x_n\}$, the power function is
$$
\begin{array}{rcl}
P^2_{X,\Phi}(x)& :=&
\displaystyle{\inf_{a\in\R^n}\left\|\delta_x-\sum_{j=1}^na_j\delta_{x_j}\right\|^2_\Phi   }. \\
\end{array} 
$$
and if we apply the scaling law \eref{eqscaleddeltas}, then
$$
P_{X,\Phi_\epsilon}(x)=P_{\epsilon X,\Phi}(\epsilon x) \fa x\in \R^d.
$$
For Power Functions, scaling points and scaling kernels is the same,
and the invariance relation is
$$
P_{X,\Phi}(x)=P_{\epsilon X,\Phi_{1/\epsilon}}(\epsilon x).
$$
\biglf
The Power Function
$P_{X,\Phi_\epsilon}(x)=P_{\epsilon X,\Phi}(\epsilon x)$
can be bounded via the standard arguments in \cite{wendland:2005-1}.
If the usual geometric properties on the domain $\Omega$ and the point set $X$
are satisfied, there is a $\Phi$- and geometry-dependent function $F_\Phi$
such that
$$
P_{X,\Phi}(x)\leq F_\Phi(h) \fa x\in \Omega
$$
for the {\em fill distance} 
$$
h:=\displaystyle{  \sup_{y\in \Omega}\min_{x \in X}\|y-x\|_2}.
$$
Scaling $X$ and $x$ to $\epsilon X$ and $\epsilon x$
leaves the geometry
invariant, and therefore 
$$
P_{\epsilon X,\Phi}(\epsilon x)\leq F_\Phi(\epsilon h)
$$
for $\epsilon\leq 1$.
The function $F_\Phi(h)$ can have algebraic or exponential decay towards
zero for $h \to 0$. A similar argument is in
\RScite{davydov-schaback:2019-1} for convergence of stencils,
though restricted to Sobolev spaces. In the case $W_2^m(\Omega)$
for $\Omega\subset\R^d$, the above function is
$F_\Phi(h)=ch^{m-d/2}$, and for kernels
satisfying Theorem\RSref{TheadmKerScaFin} it is $F_\Phi(h)=ch^{\beta/2-d/2}$.
\biglf
Note that \RScite{schaback:1995-1}
provides a two-sided bound with the same $F_\Phi$ in case of asymptotically
uniformly distributed points and kernels with finite smoothness.
%****************************************************************
\subsection{Scaling of Lagrangians}\RSlabel{SecSoL}
For Lagrange basis functions $u_j^{X,\Phi}$ we get
$$
\begin{array}{rcl}
  \delta_{j,k}=u_j^{X,\Phi_\epsilon}(x_k)&=&
 \displaystyle{   \sum_{x_i\in X} a_{ji}(X, \Phi_\epsilon)\Phi(\epsilon
   x_k-\epsilon x_i)},\\
u_j^{X,\Phi_\epsilon}(x)&=&
 \displaystyle{   \sum_{x_i\in X} a_{ji}(\epsilon X, \Phi)\Phi(\epsilon
   x-\epsilon x_i)},\\
 &=& u_j^{\epsilon X,\Phi}(\epsilon x).
\end{array}
$$
This proves that scaling points or kernels is the same thing
for calculating Lagrangians, and there is an invariance relation
$$
u_j^{X,\Phi}(x)=u_j^{\epsilon X,\Phi_{1/\epsilon}}(\epsilon x)
\hbox{    or    } u_j^{X,\Phi_\epsilon}(x)=u_j^{\epsilon X,\Phi}(\epsilon x).
$$
The paper \RScite{demarchi-schaback:2010-1} proves that Lagrangians
for kernels with limited smoothness are
uniformly bounded if points are asymptotically uniformly distributed.
Checking the proof under consideration of scaling reveals that it works
independent of scaling, because the scale factor cancels out.
The geometric assumptions are satisfied, because the form
$u_j^{\epsilon X,\Phi}(\epsilon x)$ just scales the geometry by $\epsilon$.
But the proof uses a bump function argument, and therefore fails
for analytic kernels like Gaussians or inverse multiquadrics.
This is in line with the fact that Lagrangians diverge in the flat limit
for analytic kernels if they do not converge to polynomials
\RScite{fornberg-et-al:2004-1} in degenerate cases, and it is also in
line with the convergence towards polyharmonic interpolants
\RScite{song-et-al:2012} for non-analytic kernels.
%****************************************************************
\subsection{$L_\infty$ Error Under Scaling}\RSlabel{SecLEuS}
If a user just has a single interpolation problem for a fixed point set
$X$ and wants to use the best possible kernel and scaling,
the error
$$
f(x)-s_{f,X,\Phi_\epsilon}(x)=f(x)-
\displaystyle{\sum_{x_j\in X}f(x_j)u_{x_j}^{X,\Phi_\epsilon}    } 
$$
just depends on the behavior of the Lagrangians $u_{x_j}^{X,\Phi_\epsilon}$.
For large $\epsilon$, these converge to delta functions $\delta_{x_j}(x)$
for all analytic kernels, and for $\epsilon \to 0$ they will converge
to a limit Lagrangian, except for degenerate situations in case of
analytic kernels, where they converge to $\pm \infty$. The limit Lagrangians
are polynomials for analytic kernels and polyharmonics for
kernels with finite smoothness.
\biglf
When evaluating the error as a function of
$\epsilon$, it may be that the minimum is at $\epsilon=0$,
i.e. the flat limit. In particular, this is to be expected when
$f$ is close to a polynomial and analytic kernels are used,
or if $f$ is close to a polyharmonic otherwise.
But plenty of easy experiments show that the error above may have a minimal
$L_\infty$ norm for some positive $\epsilon$, in a specific situation
determined by $f,\;X$, and $\Phi$. Observing such cases
does not exhibit a clear pattern because of the strong dependence on
$f$, and it is hard to attribute to
one of the three ingredients.
\biglf
One can get rid of the influence of $f$ by focusing on the Lagrangians,
but this will not be what single-case users want. They might not even be
interested at all in function spaces or smoothness properties of $f$,
using only the values $f(X)$. To go on from here, additional properties
of Lagrangians are needed, e.g. moment conditions, but this
leads to open problems.
%****************************************************************
\section{Scaling of Error Bounds}\RSlabel{SecSoEB}
The standard bound
\bql{eqerrstand}
|f(x)-s_{f,X,\Phi}(x)|\leq P_{X,\Phi}(x)\|f\|_\Phi
\eq
can be scaled in different ways. Scaling only $\Phi$ results in
\bql{eqphiscaleonly}
\displaystyle{\left|f(x)-  \sum_{x_j\in X}f(x_j)u_j^{X,\Phi_\epsilon}(x)   \right|}
\leq P_{X,\Phi_\epsilon}(x)\|f\|_{\Phi_\epsilon},
\eq 
playing the left-hand side
back to the Lagrangians and Section \RSlabel{SecLEuS}.
To see whether the bound gives more information,
both factors in the right-hand side need further consideration. 
Scaling only $f$ and scaling only the points are not very useful
and therefore omitted here.
\biglf
The inequality is tight for some $\epsilon$-dependent function,
but it is not clear how the ingredients vary with $\epsilon$
for a fixed $f$ and whether the two sides stay close. Recall that
Theorem \RSref{TheadmKerScaInf} excludes scales $\epsilon<1$
for arbitrary $f$ in the case of analytic kernels.
\biglf
The factors in the right-hand side were dealt with in sections
\RSref{SecSoKbN} and \RSref{SecStPF}, and we only have to look at the product.
The error bound could be improved if the native space norm
  were taken in the domain-dependent native space. But this localized norm
  is bounded above by the global one \RScite{wendland:2005-1},
  and this is the way
  the error bound is used in practice. We keep the global norm here.
\biglf
\begin{Theorem}\RSlabel{TheSobBound}
  For kernels satisfying Theorem \RSref{TheadmKerScaFin}
  and for asymptotically uniformly
placed points, the right-hand side of the standard error bound
\eref{eqerrstand} behaves like
$$
\begin{array}{rcl}
{\cal O}(1) & \hbox{ for } \epsilon\to 0,\\ 
{\cal O}(\epsilon^{\beta}) & \hbox{ for } \epsilon\to \infty.
\end{array}
$$
In Sobolev spaces $W_2^m(\R^d)$ with integer $m$,
the square of the bound behaves like
$$
\displaystyle{ 
\epsilon^{2m}\sum_{j=0}^m{m\choose j}\epsilon^{-2j}
|f|^2_{W_2^j(\R^d)}}  \hbox{ for } \epsilon\to 0.
$$\qed
\end{Theorem} 
\begin{proof}
From the previous section, we get that the error bound is ${\cal O}(1)$ for $\epsilon \to 0$,
caring for  the flat limit. For large $\epsilon$, the behavior
of both factors leads
to the ${\cal O}(\epsilon^{\beta})$ asymptotics of the bound.
\end{proof}
So far, this result is disappointing, because it is not clear
if there is a minimum with respect to $\epsilon$.
Some experimental observations are in Section
\RSref{SecEBaFoS} showing typical cases in Figure \ref{FigscAEB441}.
%****************************************************************
\section{Natural Scales of Functions}\RSlabel{SecNSoF}
To proceed towards finding a ``natural'' scale of a global function
$f \in \calh_\Phi$ with respect to the kernel chosen,
we consider the minimum of $\|f\|_{\Phi_\epsilon}$ over all scales $\epsilon $
that are $\Phi$-admissible for $f$. By \eref{eqfPlowbnd}, the search for optimal scales can be limited to
\bql{eqepslimit}
\epsilon^d\leq\dfrac{\|f\|^2_{\Phi}\|\hat\Phi\|_\infty}{\|f\|^2_2}.
\eq
For kernels with finite smoothness, Theorem \RSref{TheadmKerScaFin}
guarantees such a minimum, because $\|f\|_{\Phi_\epsilon}$
goes to infinity like $\epsilon^d$ for $\epsilon\to\infty$,
and like $\epsilon^{d-\beta}$ for $\epsilon\to 0$ with $\beta >d$.

For general functions and analytic kernels, Theorem
\RSref{TheadmKerScaInf} limits scales to $\epsilon\geq 1$,
while \eref{eqepslimit} gives an upper limit in all cases.
  Therefore there must be a minimum of $\|f\|_{\Phi_\epsilon}$
  in this restricted domain.
But things are strongly $f$-dependent,
and if $f$ is bandlimited,  we can look at the limit $\epsilon\to 0$ even
  for analytic kernels. Then 
  \eref{eqfeffe} 
  shows that the left bound gets constant for
  $\epsilon\to 0$, but the right one
  will blow up. There may be a minimum at
  $\epsilon=0$, i.e. for the flat limit case. 
\biglf
By the scaling law $\|f_\epsilon\|_\Phi^2=\|f\|^2_{\Phi_{1/\epsilon}}$,
the minimization of $\|f_\epsilon\|_\Phi^2$ is equivalent to minimization
of $\|f\|^2_{\Phi_{1/\epsilon}}$, and is therefore also covered by the above
argumentation, except that $\epsilon$ must be replaced
by $1/\epsilon$.
\biglf
Summarizing, all functions in native spaces of
a kernel with finite smoothness have a natural scale
as an $\epsilon$ that minimizes $\|f\|^2_{\Phi_{\epsilon}}$.
For analytic kernels, existence of a minimum 
of $\|f\|_{\Phi_\epsilon}$ will be crucially dependent on $f$, but it may arise
in the flat limit $\epsilon=0$.
Even if it could be estimated efficiently and accurately, it
would not necessarily lead to a minimal error in
\eref{eqphiscaleonly} for specific cases. But it may have a
small factor in front of the convergence rate
associated with the unscaled kernel. Some experimental results are in Section
\RSref{SecNSNaFoS}, but they replace $\|f\|_{\Phi_\epsilon}$ by
$\|s_{f,X,\Phi_\epsilon}\|_{\Phi_\epsilon}$ and the resulting scales are not
particularly useful. 
%****************************************************************
\section{Examples}\RSlabel{SecEx}
A comparison of existing scale estimation algorithms
is not attempted here. The examples will be
strictly confined to illustrations of
the theoretical arguments of this paper. Sophisticated algorithms are avoided
in favor of what scientists do when they want quick results.
Throughout, we work in the unit square in $\R^2$, and we limit all
hazardous scales in a primitive way by stopping when MATLAB's
{\tt condest} exceeds $10^{14}$. In contrast to many other papers, we
add the flat limit situation, not as a limit, but by adding a polynomial
and a polyharmonic solution.  
%****************************************************************
\subsection{Remarks on Kernels and Functions}\RSlabel{SecRoK}
To keep the number of final figures small,
only four kernels and functions were selected as typical cases.
\biglf
From the class of analytic kernels, we only consider
the Gaussian and the inverse multiquadric
$\Phi(x,y)=(1+\|x-y\|_2)^{-1/2}$. 
For analytic functions, we only consider the scaled class
$f_\epsilon(x)=(1+\epsilon^2\|x\|_2^2)^{-1}$
that models the Runge phenomenon. Their complex extensions have
singularities at distance $1/\epsilon$ from the real axis, and we use
$\epsilon=1/4$ for a good and $\epsilon=25$ for a bad case.
Both are polynomials up to machine accuracy, but the bad case needs a
very high degree for this. One can expect that the polynomial flat limit works
fine for the good case and fails for the bad one, even though the function
is analytic on the real numbers.
\biglf
The non-analytic case has to compromise on the kernel side.
The intention is to
use kernels with comparable $W_2^3(\R^2)$ smoothness.
The matching 2D polyharmonics are
$\|x-y\|_2^{4}\log\|x-y\|_2$, and the matching
radial Mat\'ern-Sobolev kernel is $K_2(r)r^2$. But standard Wendland kernels
$\phi_{3,k} $ on $\R^2$ work in $W_2^{k+3/2}(\R^2)$ and cannot match
$W_2^3(\R^2)$. Picking $k=2$, i.e.
the radial kernel $(1-r)_+^6(35r^2+18r+3)$ leads to $W_2^{3.5}(\R^2)$
smoothness, and therefore we also use the polyharmonic
kernel $\|x-y\|_2^{6}\log\|x-y\|_2$ working in $W_2^4(\R^2)$. 
\biglf
In legends of plots, kernels are abbreviated by {\tt g},  {\tt mq},
{\tt ms3},  {\tt w3.5}, and for the flat limit we use
{\tt p}, {\tt ph3}, and {\tt ph4} for polynomials and polyharmonics.
The trailing number indicates Sobolev smoothness, not the usual RBF parameters.
To match the basic kernel scale to the domain $[-1,+1]^2$
and regular point sets of 121 or 441 points, the Gaussian, the inverse
multiquadric
and the Wendland function were pre-scaled by $\epsilon=10, 10$, and $0.2$,
respectively. Note that for small $\epsilon$ the Wendland functions
have large supports, and the matrices are not sparse. Error evaluation is
made in cell midpoints only. 
\biglf
Polynomials were fitted to the data by brute force, multiplying data
values by the pseudoinverse of a high-order Vandermonde matrix.
This is a hazardous
thing by itself, but we wanted to use primitive techniques only.
Since the result deteriorates numerically
when large polynomial degrees are taken, we picked a degree that yields a
minimal $L_\infty$ error. Of course, there are much better techniques
for polynomial interpolation.
\biglf
For non-analytic functions we picked $\|x\|_2^3$ for a case
with a single interior point of non-smoothness, and $\|x\|_\infty$
for a piecewise linear continuous case with derivative discontinuities along
lines. The chosen kernels and functions are the result of
many test runs suppressed here.
%****************************************************************
\subsection{$L_\infty$ Errors as Functions of Scales}\RSlabel{SecLEaFoS}
Figure \RSref{FigscErr121} shows observed $L_\infty$ errors
as functions of scales, for the four functions and kernels.
The point set $X$ consisted of 121 regular points on $[-1,+1]^2$,
while Figure \RSref{FigscErr441} is for 441 points, just to demonstrate
the dependence on $X$. The errors for polynomials and polyharmonics
are given as horizontal lines with markers.
\biglf
Here is a list of observations.
\begin{enumerate}
\item The top two functions are analytic, but only the left one
  is close to a low-degree polynomial. There, the flat limit for analytic
  kernels is optimal, but the non-analytic kernels do not
  perform badly, as well as the polyharmonics.
  \item The top right case shows the Runge effect. Well-chosen scales 
    of analytic kernels may give quite some improvement over the others,
    but the scale estimation will be hazardous. Non-analytic kernels
    and polyharmonics work fine when data sets get dense.
  \item The bottom left case is a mildly nonsmooth non-polynomial function
    with the non-smoothness located just in the center of the domain.
    Polyharmonics are the best choice. 
  \item The bottom right case is continuous, but nondifferentiable
    along lines. Except for lucky scales, analytic kernels should be avoided.
    The others work well.
  \item The hunt for error-minimizing scales may pay off for analytic kernels,
    but is a risky thing that often is not much better than what
    kernels with low smoothness or flat-limit polyharmonics can do.
  \item If one calculates a polynomial
    and a polyharmonic solution to pick the better one, there is
    not much to gain by bridging the gap
    between the flat limit and the smallest $\epsilon$ that works
    within a crude {\tt condest} limit. Continuing the curves to
    the left will not give something much better.
  \item Motivated by the result that small errors imply large conditions
    of kernel matrices \RScite{schaback:1995-1}, users were tempted to
    take the smallest $\epsilon$ they could handle. The right-hand plots
    show that this is misleading for analytic kernels.
    For non-analytic kernels it works fine, but the polyharmonic flat limit
    will be a better and faster solution, realizing the smallest
    $\epsilon$ ever. 
    \item
    The vast literature on algorithms for scale estimation has no clear winner
    yet. A possible reason is that the above
    plots show widely varying situations to care for. It is suggested to
    include $\epsilon=0$ into testing of scales, because 
    the flat limit solution, polynomial or polyharmonic,
    may be hard to beat by any other scaling. 
\end{enumerate} 

\begin{figure}
\begin{center}
\includegraphics[height=5.5cm, width=5.5cm]{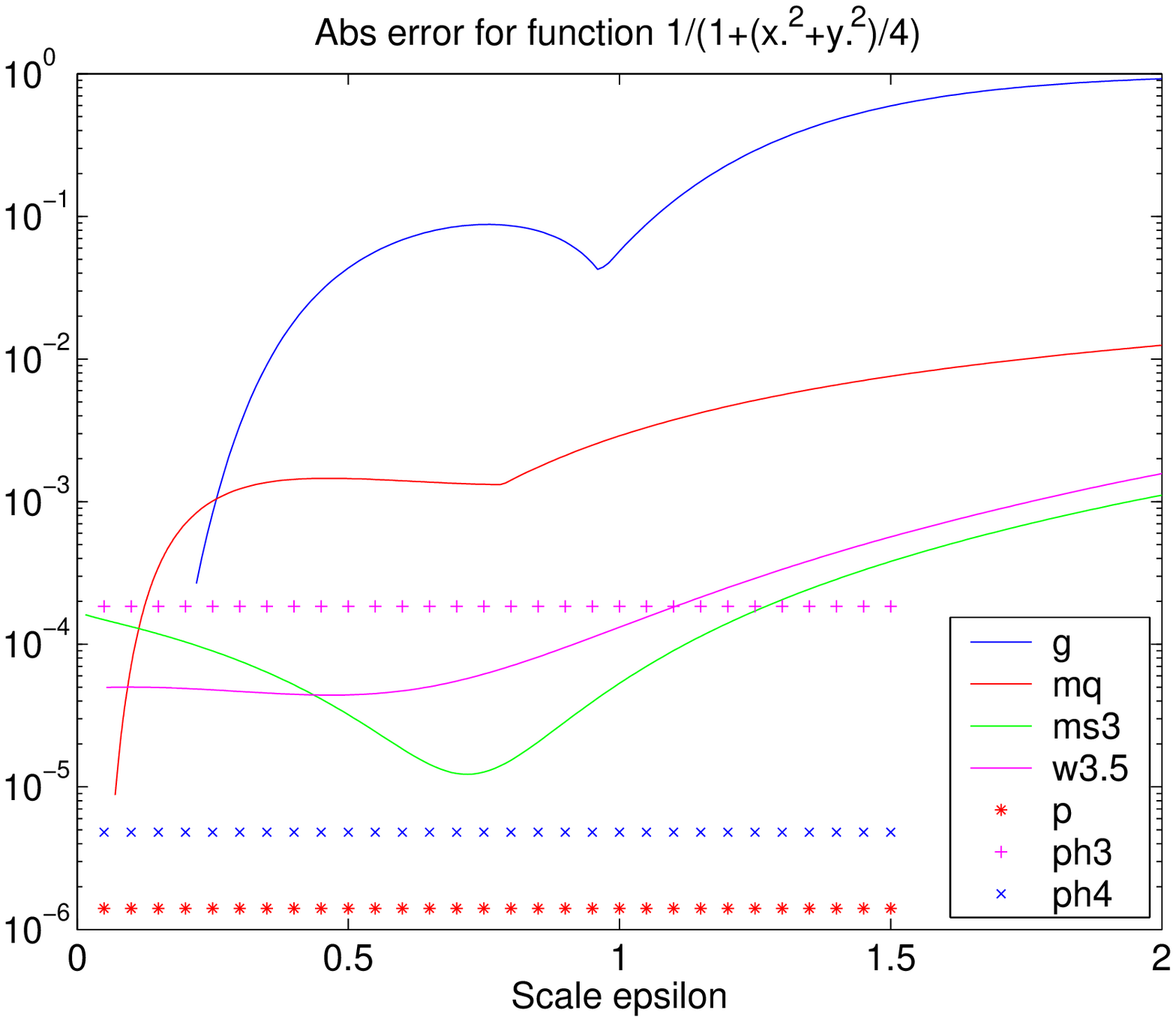}
\includegraphics[height=5.5cm, width=5.5cm]{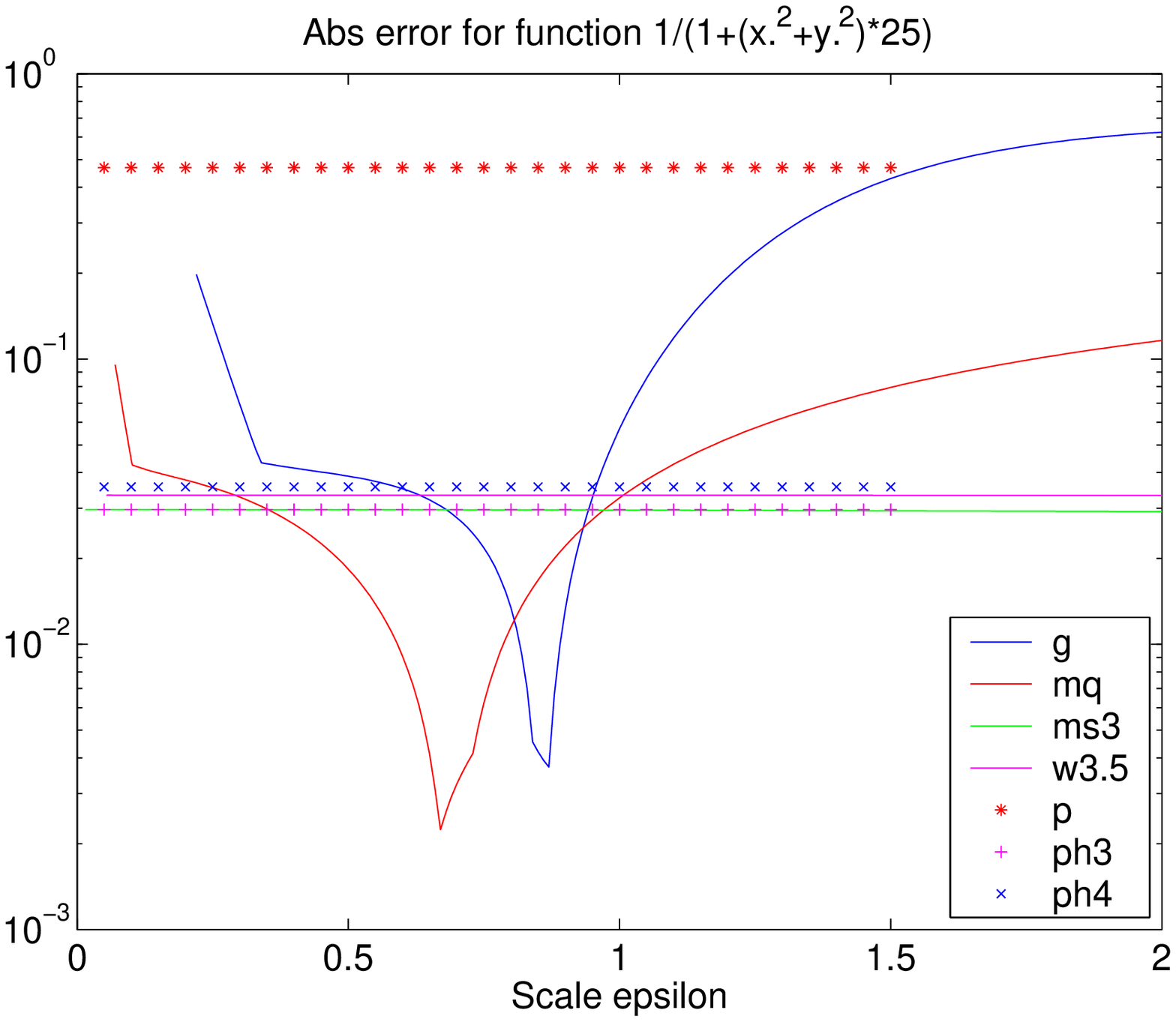}\\
\includegraphics[height=5.5cm, width=5.5cm]{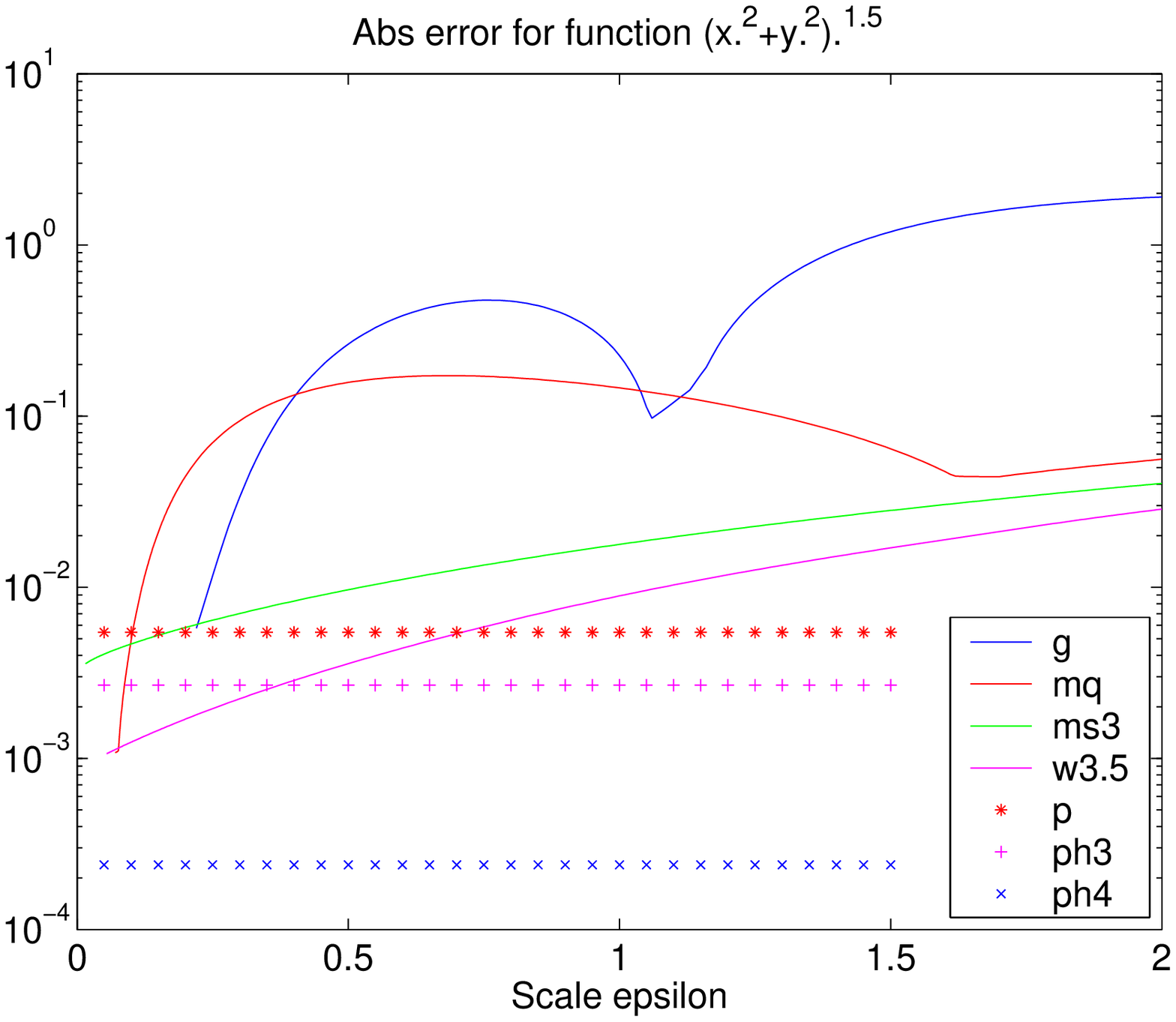}
\includegraphics[height=5.5cm, width=5.5cm]{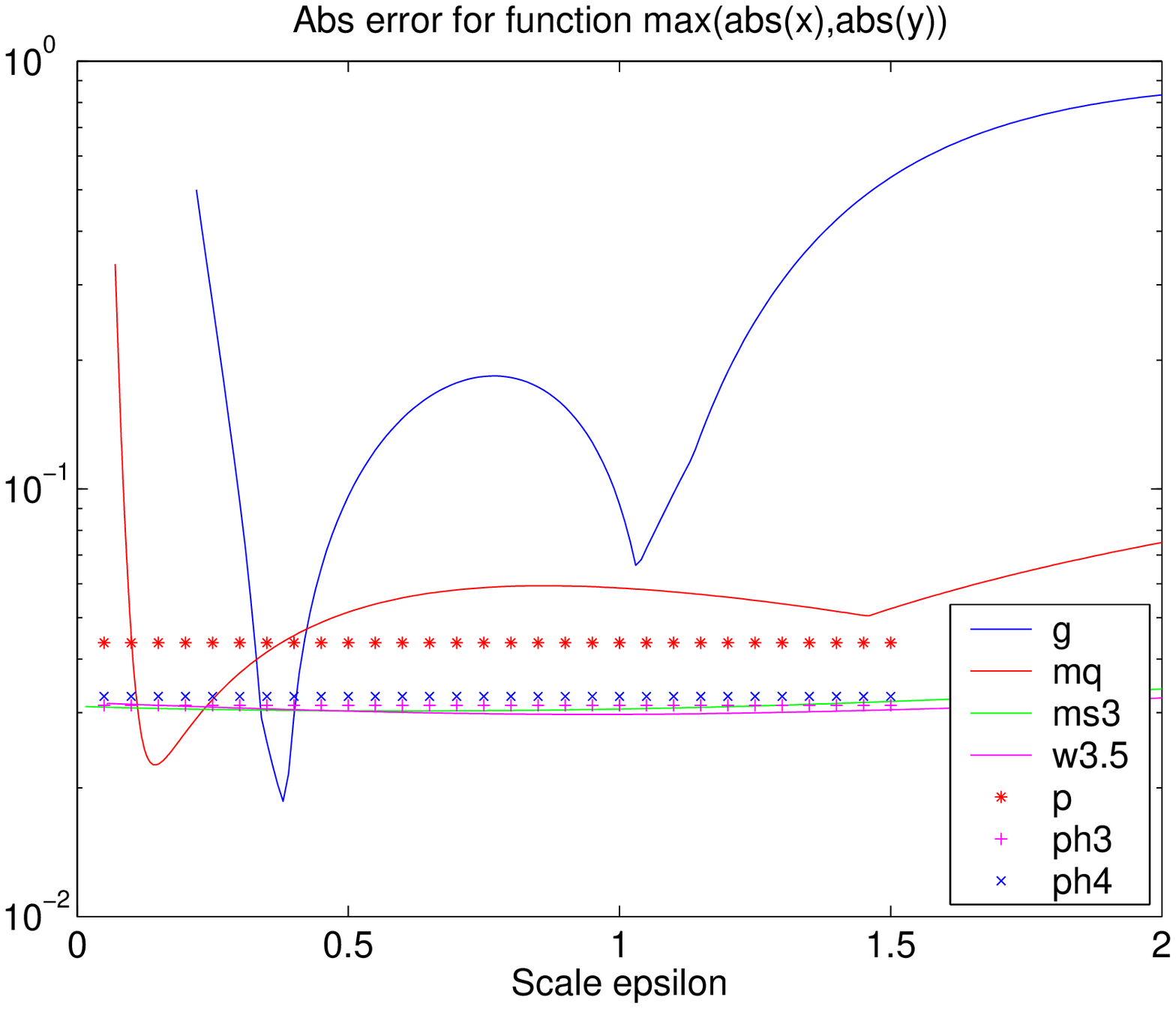}
\end{center}
\caption{$\ell_\infty$ errors on midpoints as functions of scale for
    121 regular data locations on $[-1,+1]^2$.\RSlabel{FigscErr121}}
\end{figure}

\begin{figure}
\begin{center}
\includegraphics[height=5.5cm, width=5.5cm]{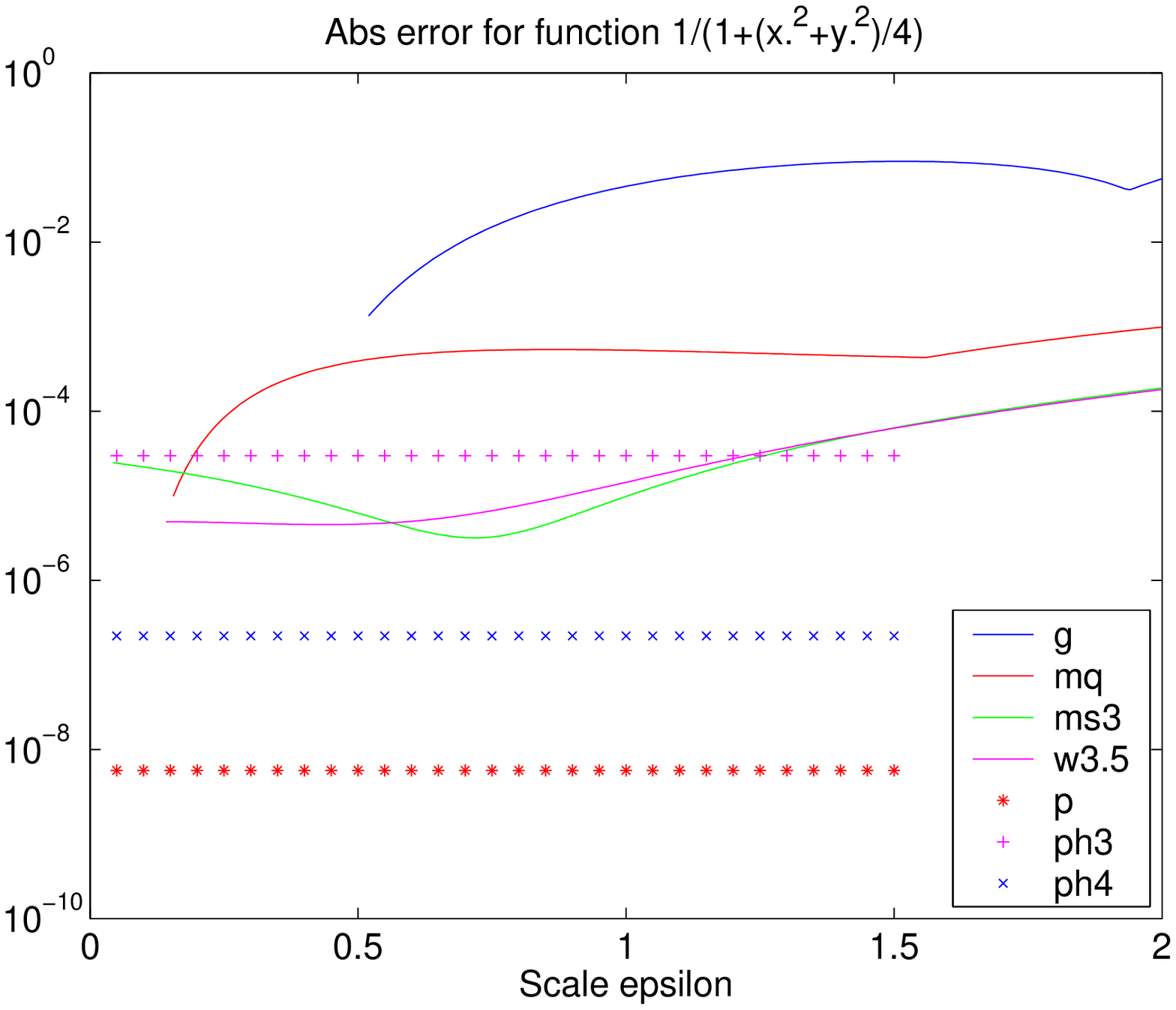}
\includegraphics[height=5.5cm, width=5.5cm]{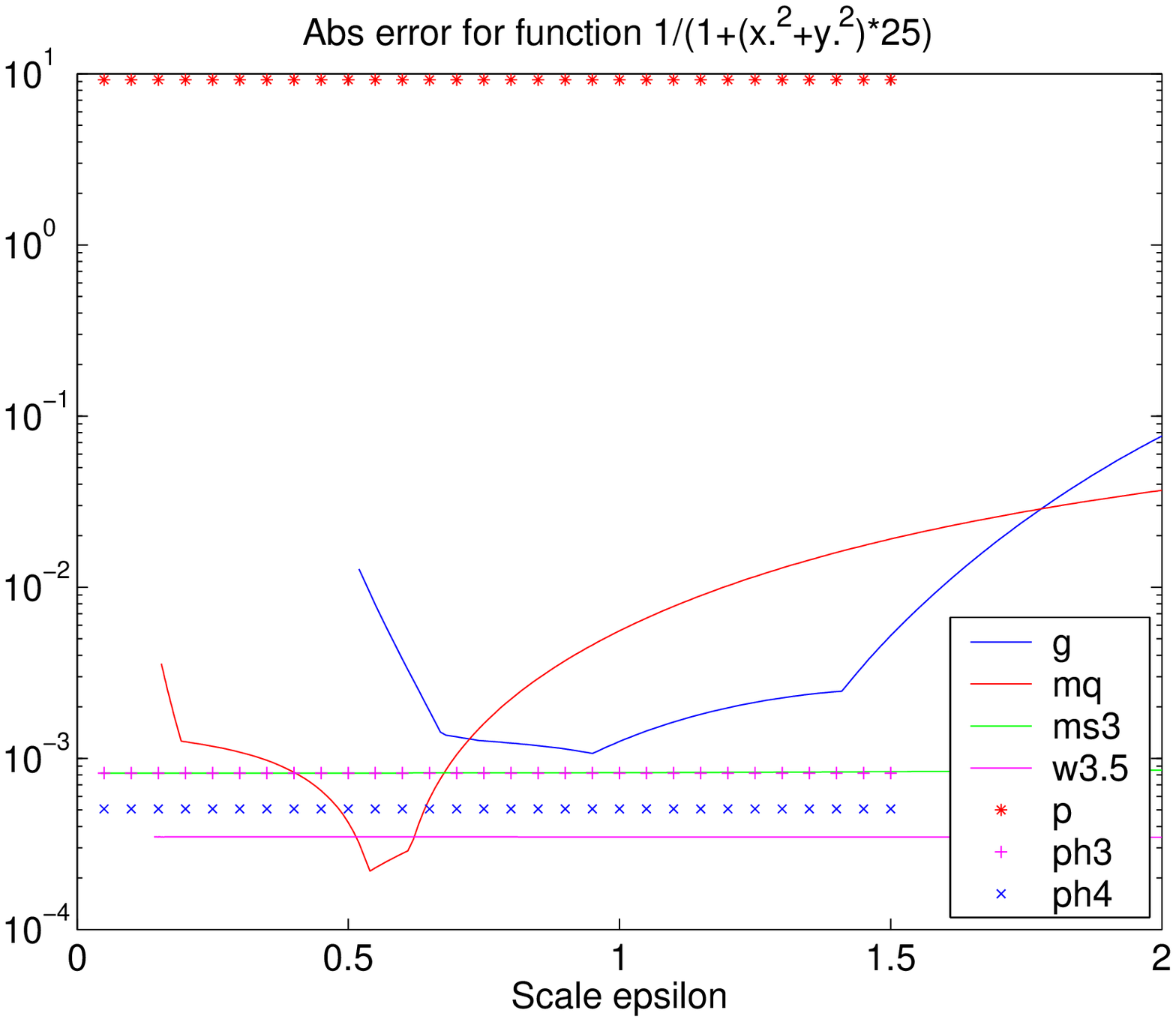}\\
\includegraphics[height=5.5cm, width=5.5cm]{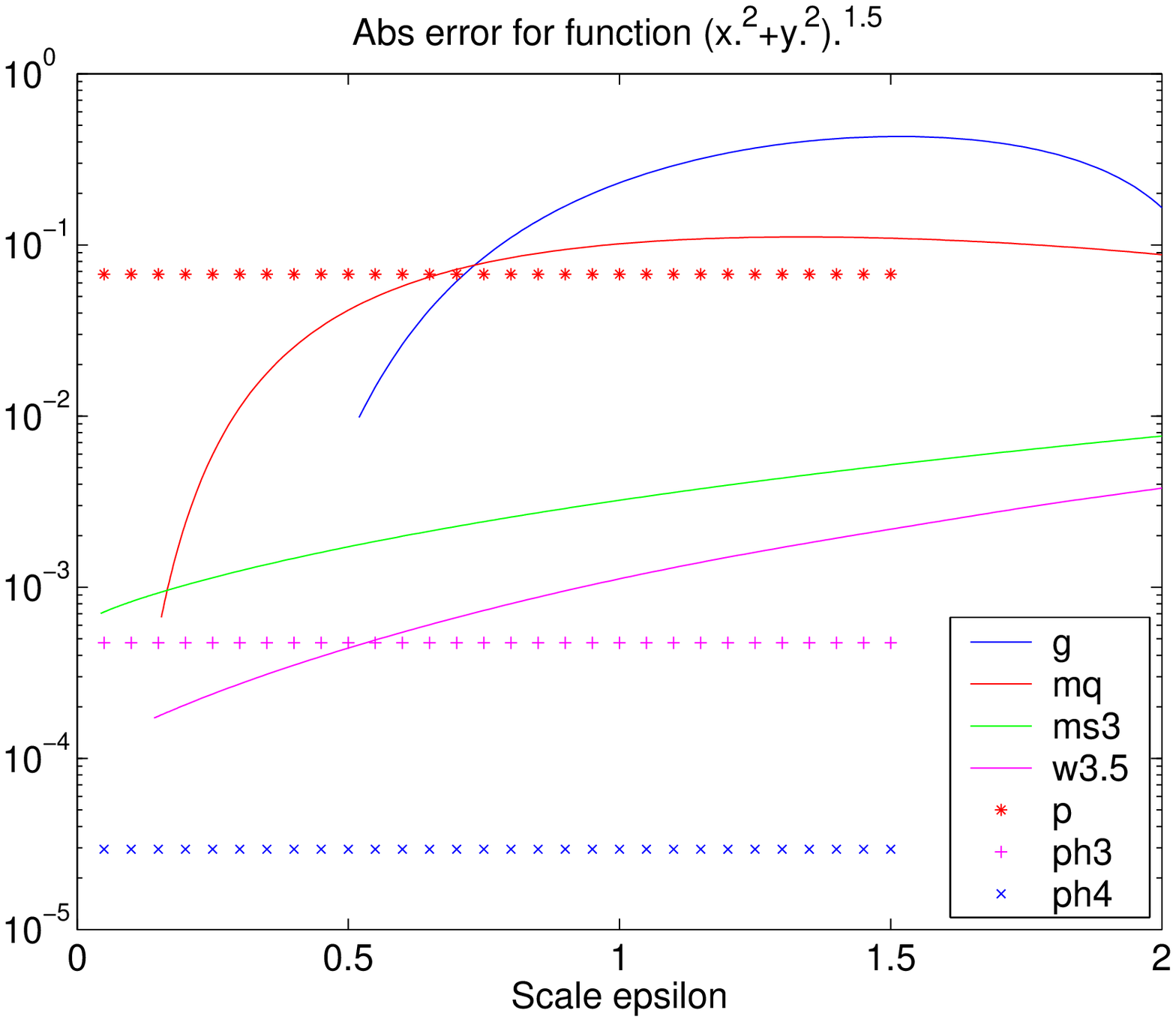}
\includegraphics[height=5.5cm, width=5.5cm]{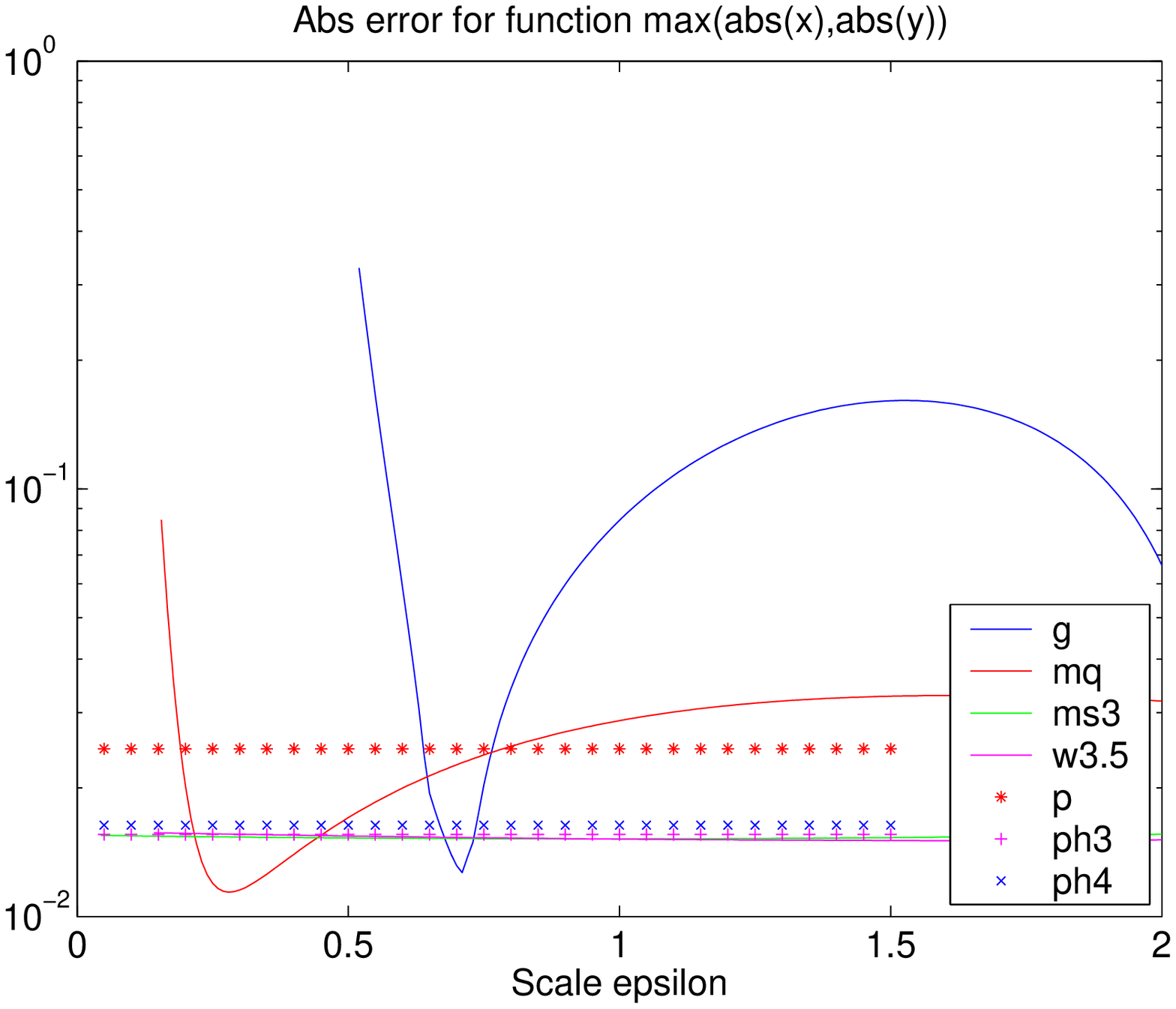}
\end{center}
\caption{$\ell_\infty$ errors on midpoints as functions of scale for
    441 regular data locations on $[-1,+1]^2$.\RSlabel{FigscErr441}}
\end{figure}
%****************************************************************
\subsection{Native Space Norms as Functions of Scales}\RSlabel{SecNSNaFoS}
Section \RSref{SecNSoF} calls for plots
of $\|f\|_{\Phi_\epsilon}$ as functions of $\epsilon$ with clear minima.
This would be another estimation method for scales, if it were
cheaply and accurately available.
But only native space norms
$\|s_{f,X,\Phi_\epsilon}\|_{\Phi_\epsilon}\leq \|f\|_{\Phi_\epsilon}$
of interpolants on point sets $X$ are available, and it is not clear
how close they are to $\|f\|_{\Phi_\epsilon}$ and whether they have the
same behavior for varying $\epsilon$ as predicted for
$\|f\|_{\Phi_\epsilon}$ by Section \RSref{SecNSoF},
in particular Theorem \RSref{TheadmKerScaFin}. There, as functions
of $c=1/\epsilon$, the quantity $c^d\|f\|_{\Phi_{1/c}}$ should be a polynomial
in $c$ of degree $m$, if $f$ is in $W_2^m(\R^d)$ and $\Phi$
is the kernel for that space. At least, one should look at
$c^d\|s_{f,X,\Phi_{1/c}}\|_{\Phi_{1/c}}$ as a function of
large $c$ for various kernels. This is still open, and maybe it
yields an estimate for a good smoothness parameter. 
\biglf
Figure \RSref{FigscNsn441} shows $\|s_{f,X,\Phi_\epsilon}\|_{\Phi_\epsilon}$
as a function of $\epsilon$ for the four selected cases in the 441 point
situation. The scale-independent native space norms of the polyharmonics
  are given for comparison as horizontal lines with markers.
There is no transition to the flat limit, because
$\|s_{f,X,\Phi_\epsilon}\|_{\Phi_\epsilon}$ seems to behave similarly to
$\|f\|_{\Phi_\epsilon}$, i.e. going to infinity in the flat limit,
exponentially for analytic kernels and algebraically for the others.
Minima of the native space norm, if there are any,
do not appear to be useful as criteria for good scales. 
\begin{figure}
\begin{center}
\includegraphics[height=5.5cm, width=5.5cm]{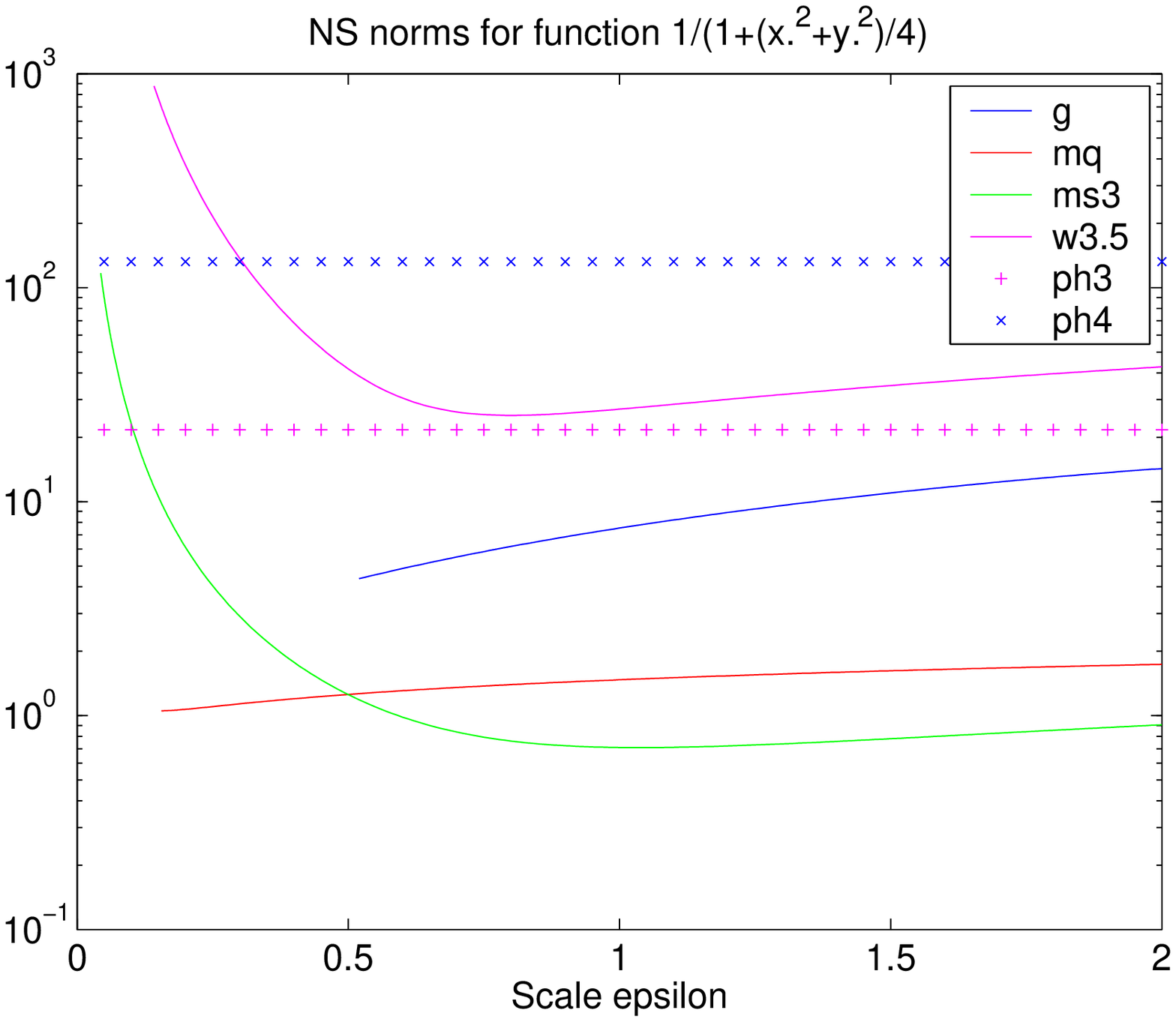}
\includegraphics[height=5.5cm, width=5.5cm]{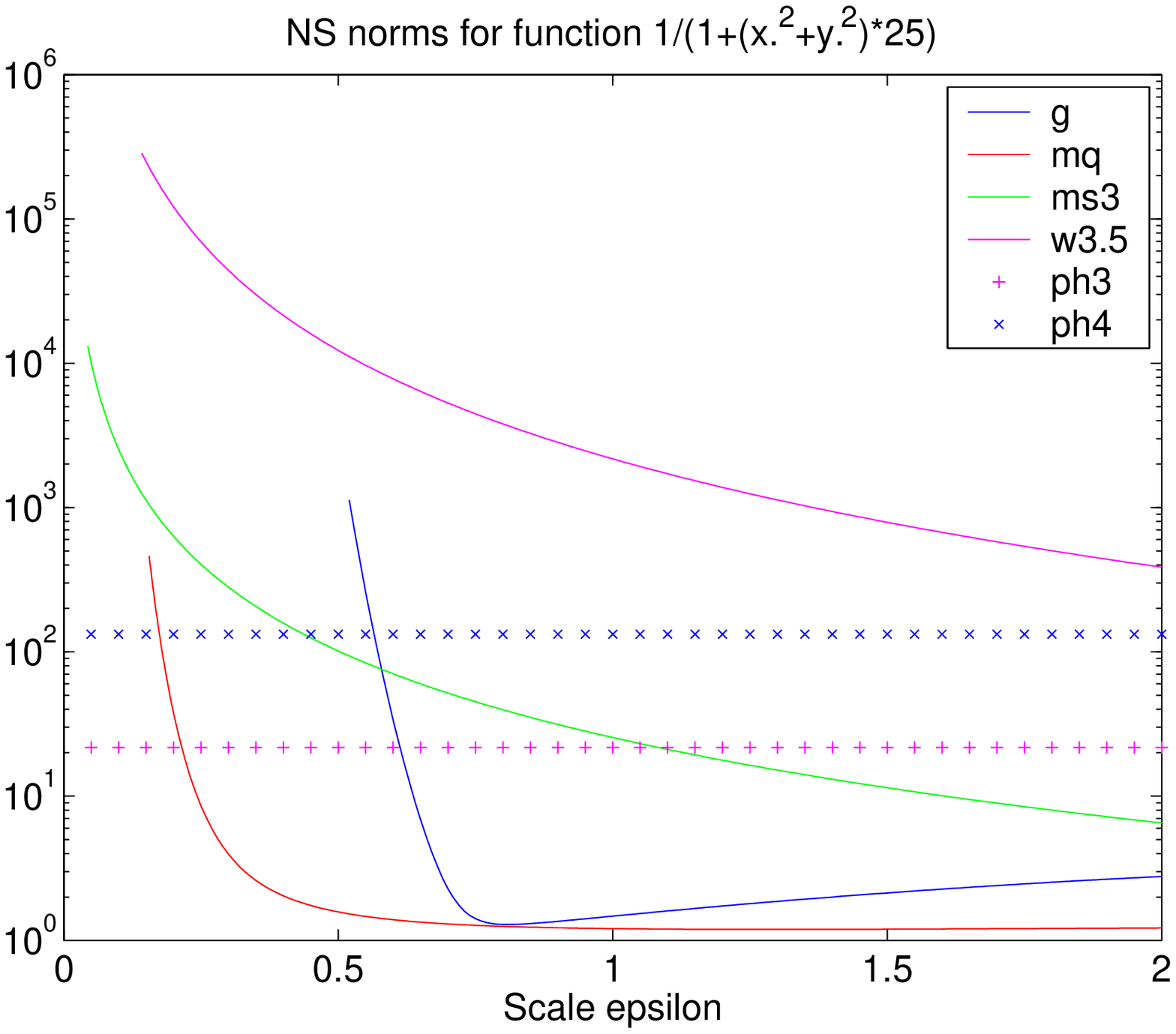}\\
\includegraphics[height=5.5cm, width=5.5cm]{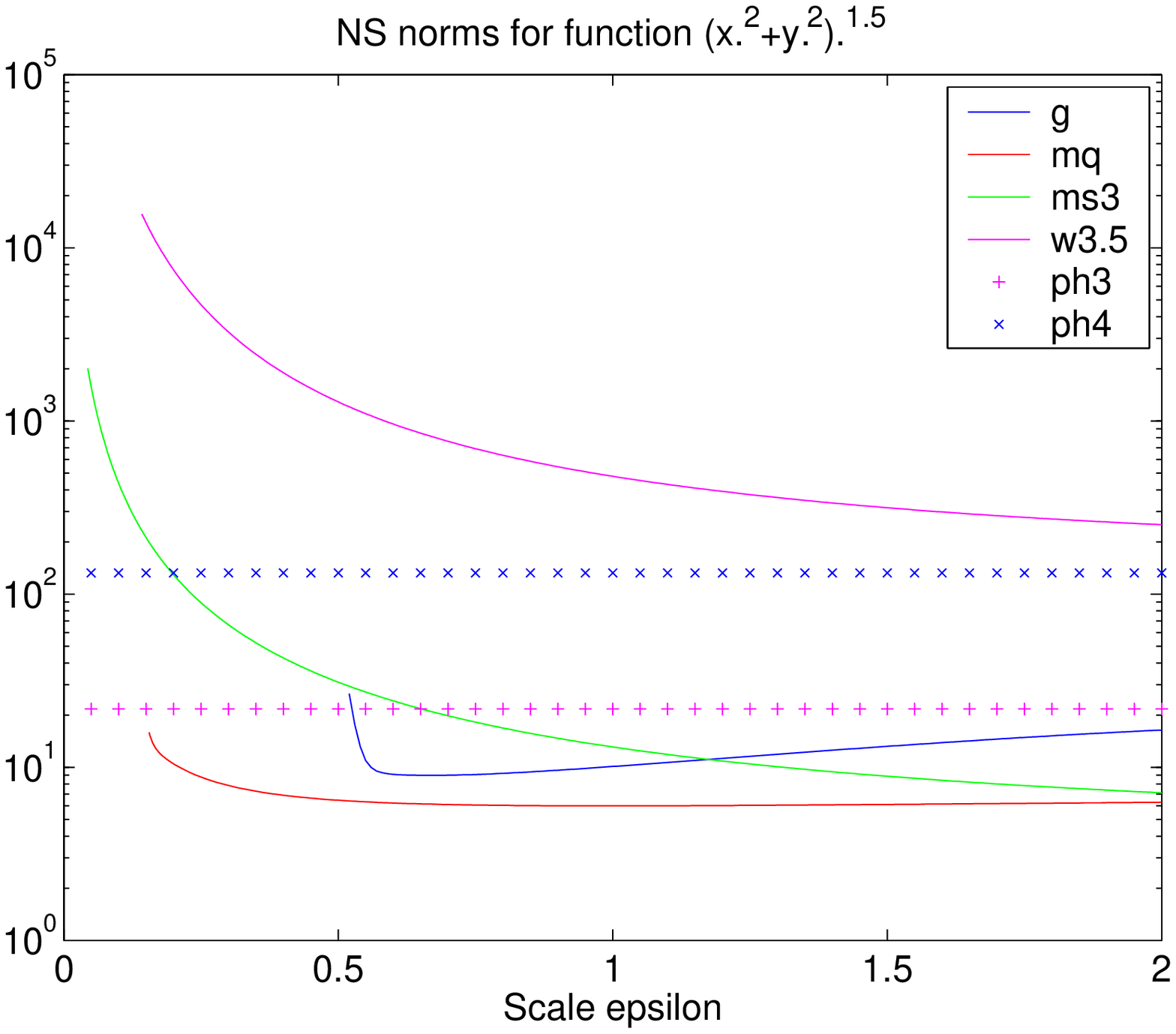}
\includegraphics[height=5.5cm, width=5.5cm]{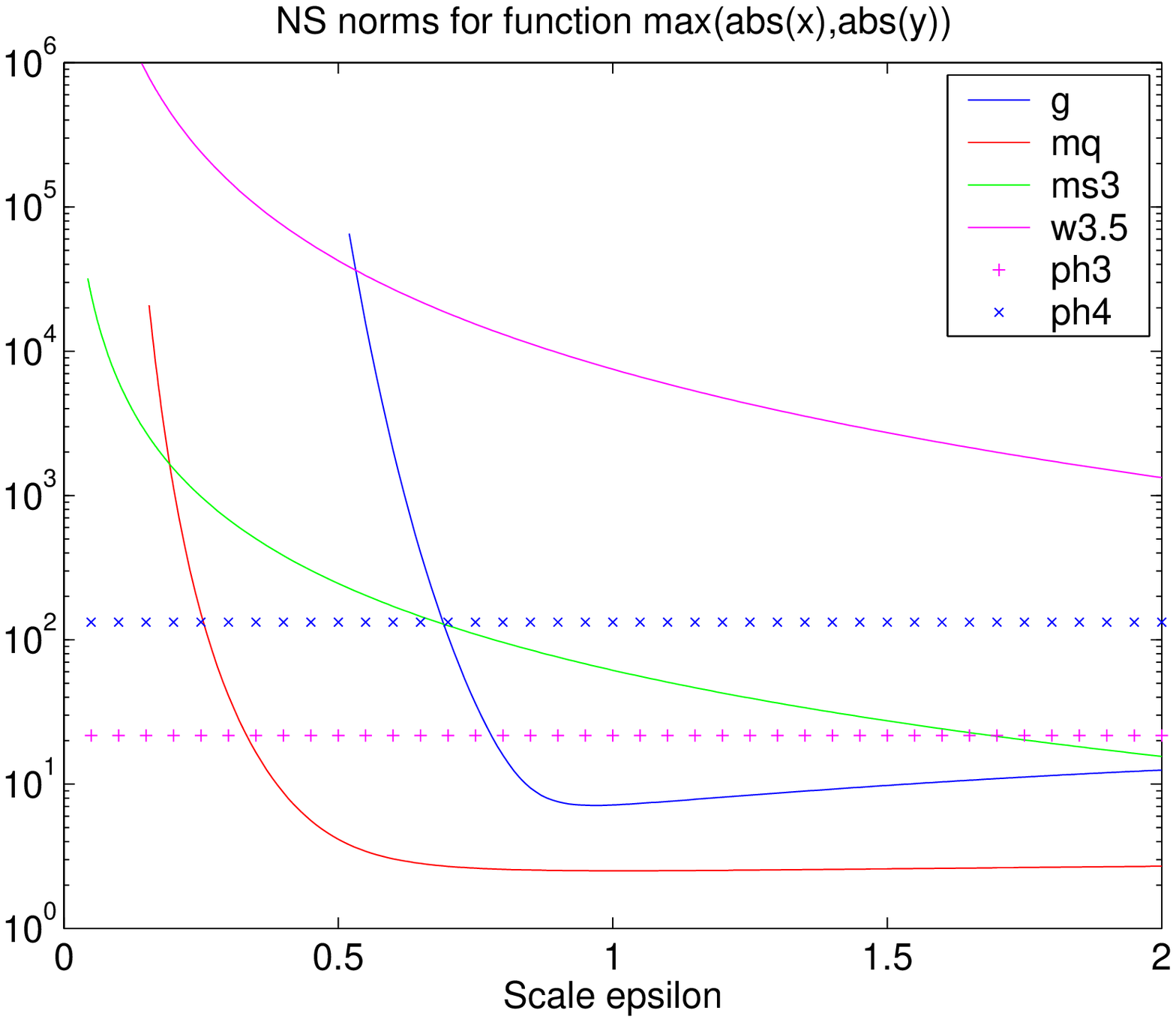}
\end{center}
\caption{Native space norms for varying scale.\RSlabel{FigscNsn441}}
\end{figure}
%****************************************************************
\subsection{Error Bounds as Functions of Scales}\RSlabel{SecEBaFoS}
Figure \RSref{FigscAEB441} shows the accessible part of the right-hand side
of the standard error bound \RSref{eqfsPf}. The replacement of
$\|f\|_{\Phi_\epsilon}$ by $\|s_{f,X,\Phi_\epsilon}\|_{\Phi_\epsilon}$ deletes
the reliability as an upper bound,  but it is interesting to see
whether it indicates good scales when compared to Figure \RSref{FigscErr441}
for the error. Note that the curves in Figure \RSref{FigscAEB441} can be
calculated without having any additional points for error evaluation.
The qualitative behavior is unexpectedly similar, and the observations made
after seeing the error plots could have been made without evaluating errors.
\begin{figure}
\begin{center}
\includegraphics[height=5.5cm, width=5.5cm]{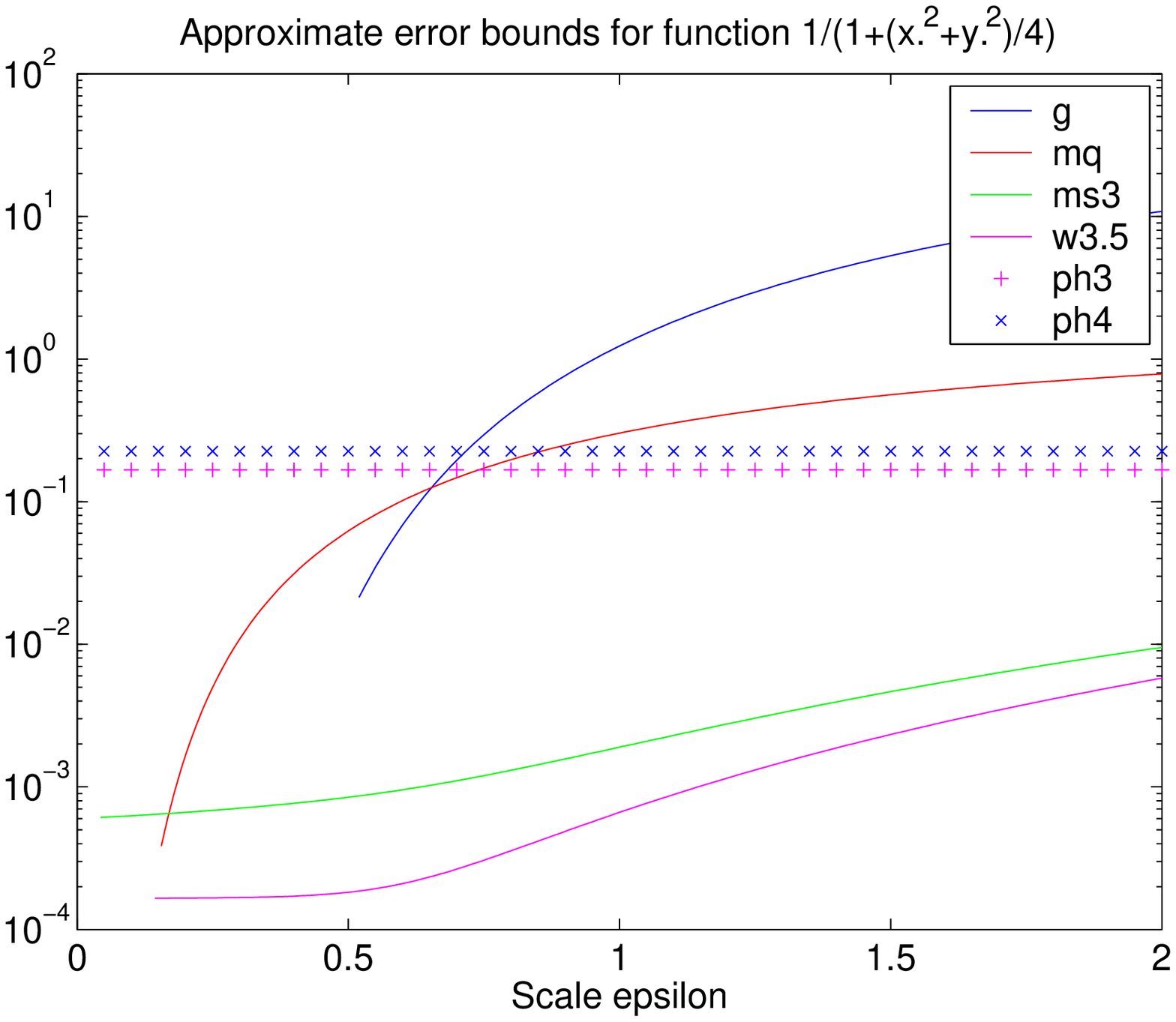}
\includegraphics[height=5.5cm, width=5.5cm]{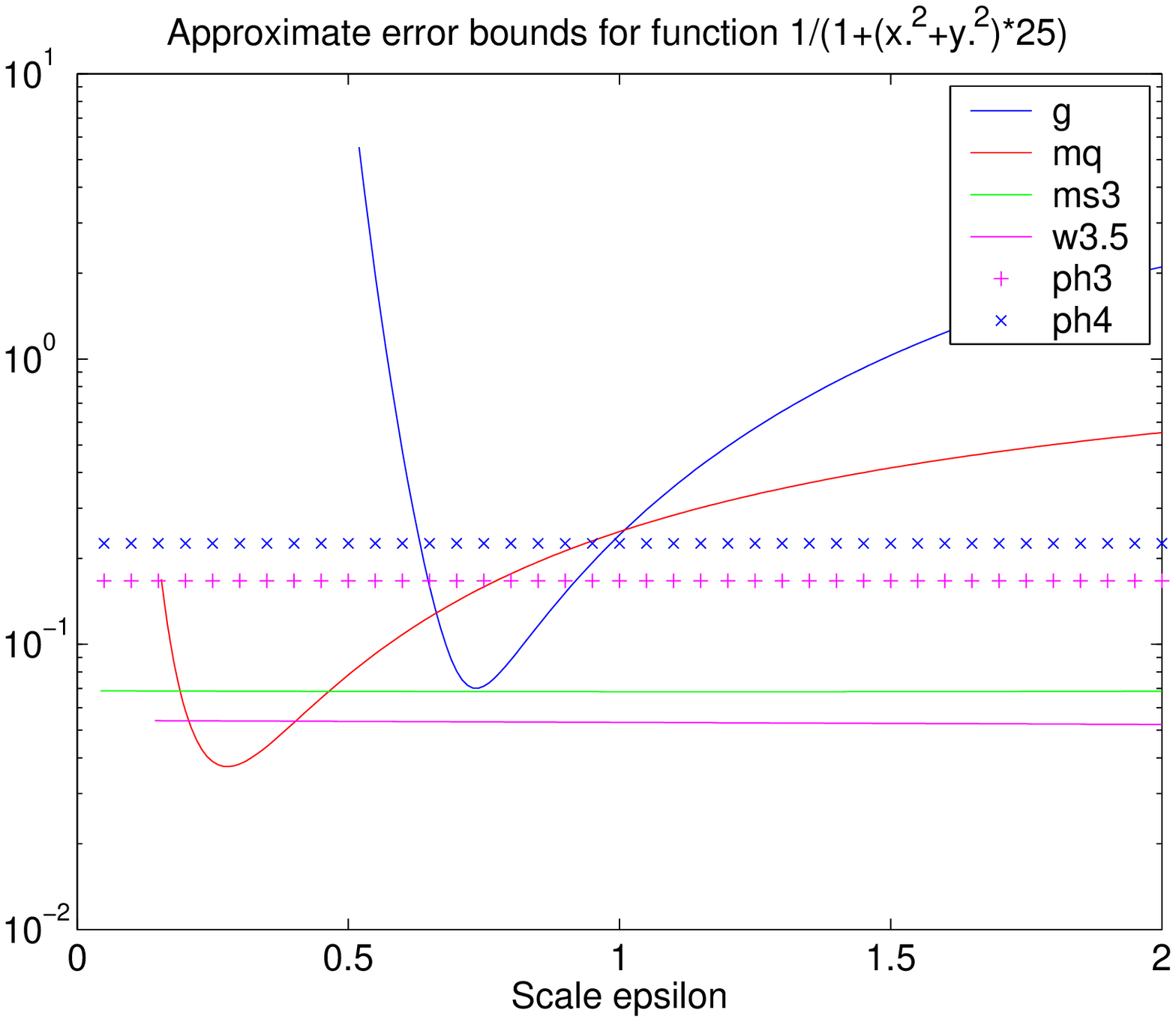}\\
\includegraphics[height=5.5cm, width=5.5cm]{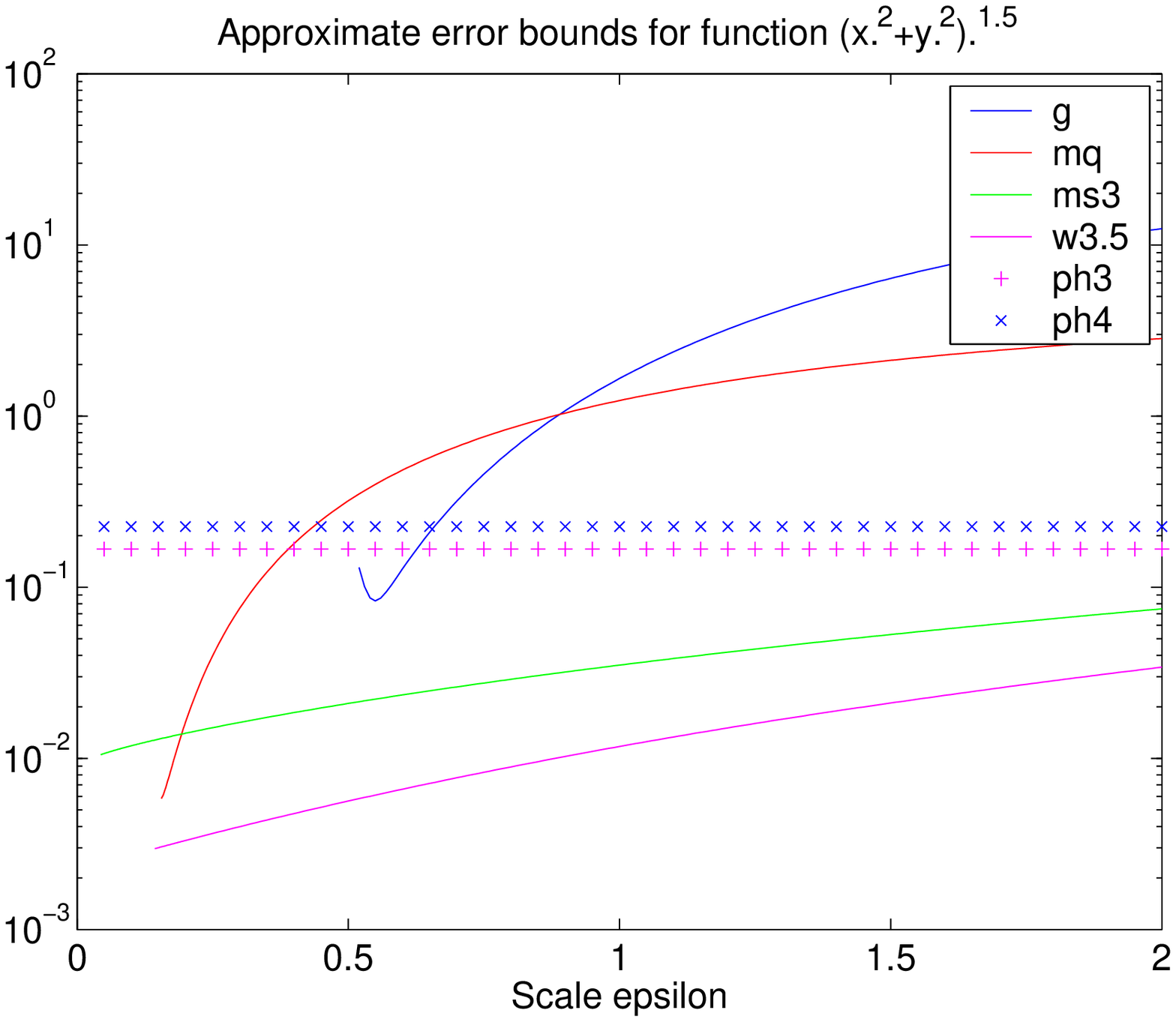}
\includegraphics[height=5.5cm, width=5.5cm]{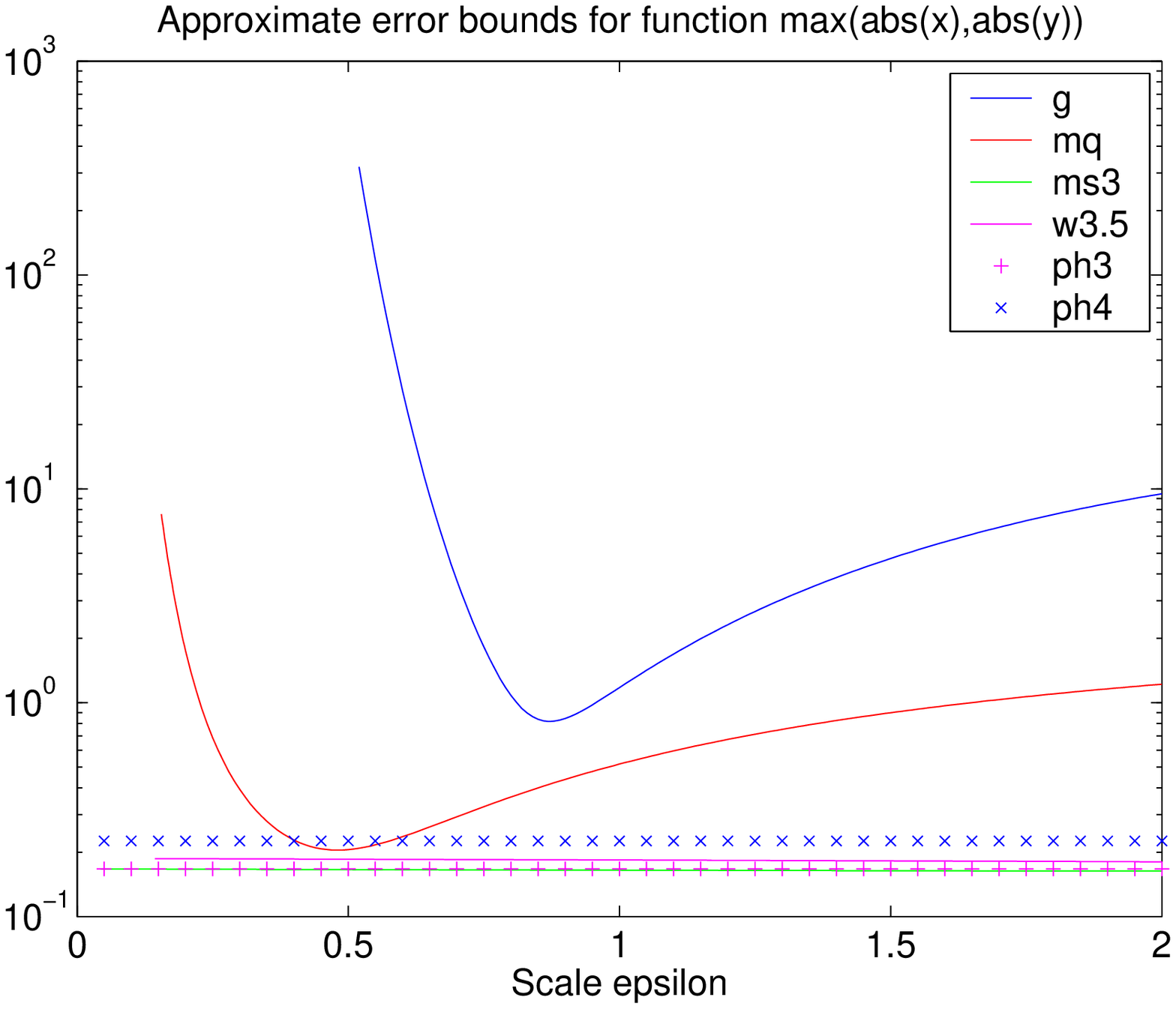}
\end{center}
\caption{Approximate error bounds for varying scale.\RSlabel{FigscAEB441}}
\end{figure}
\biglf
Figures \RSref{Fig76allinone} and \RSref{Fig28allinone} combine Figures
\RSref{FigscErr441} and  \RSref{FigscAEB441} for easier comparison
of actual errors (top row) to native space norms (center row)
and error bounds (bottom row). Figure  \RSref{Fig28allinone} has two
moderately 
well-behaving cases, while Figure \RSref{Fig76allinone} shows
results for the nonsmooth function 
$\max(|x|,|y|)$ to the right and
$1/(1+25(x^2+y^2))$ with an expectable Runge phenomenon on the left. 
\begin{figure}
\begin{center}
\includegraphics[height=5.5cm, width=5.5cm]{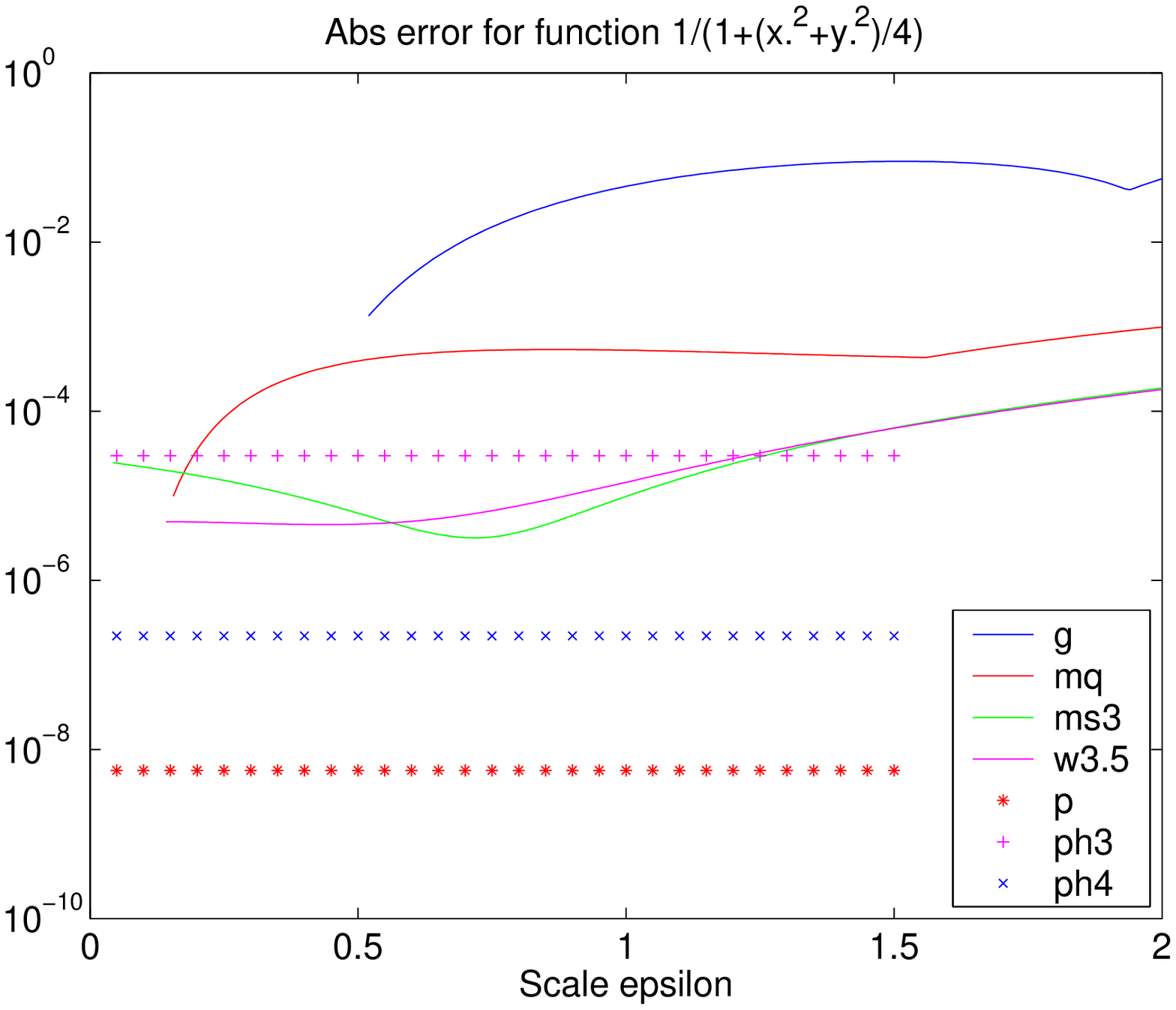}
\includegraphics[height=5.5cm, width=5.5cm]{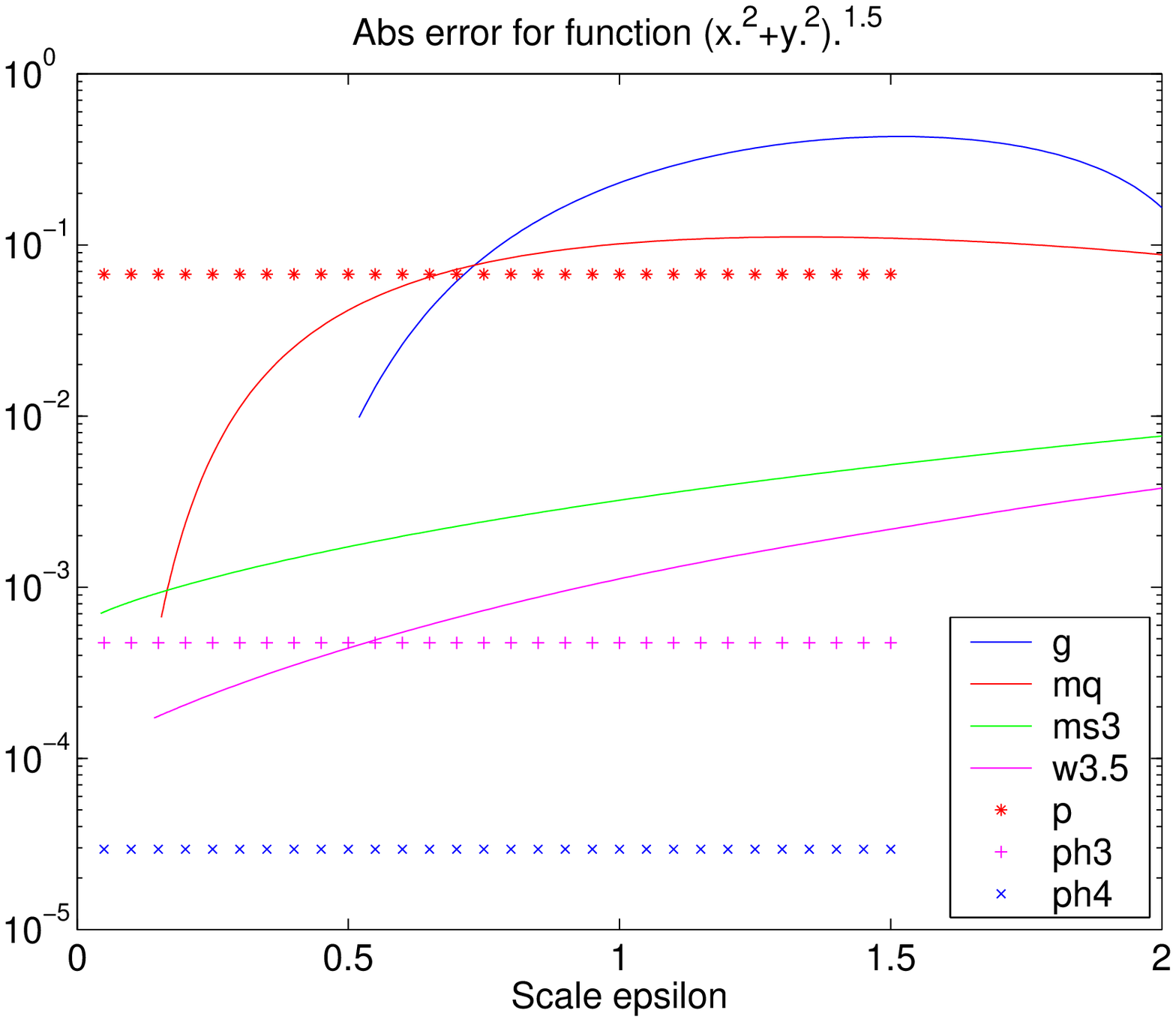}\\
\includegraphics[height=5.5cm, width=5.5cm]{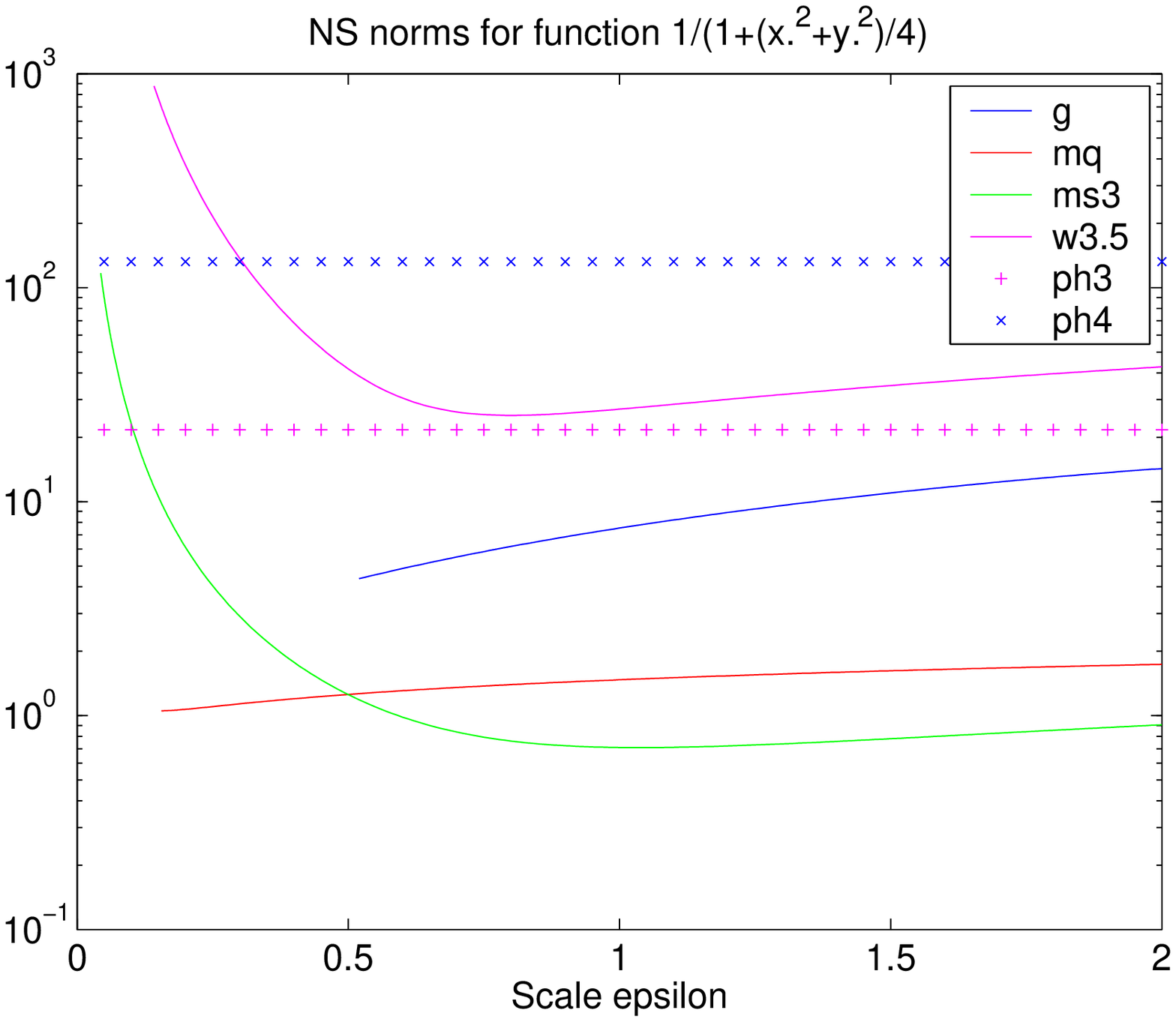}
\includegraphics[height=5.5cm, width=5.5cm]{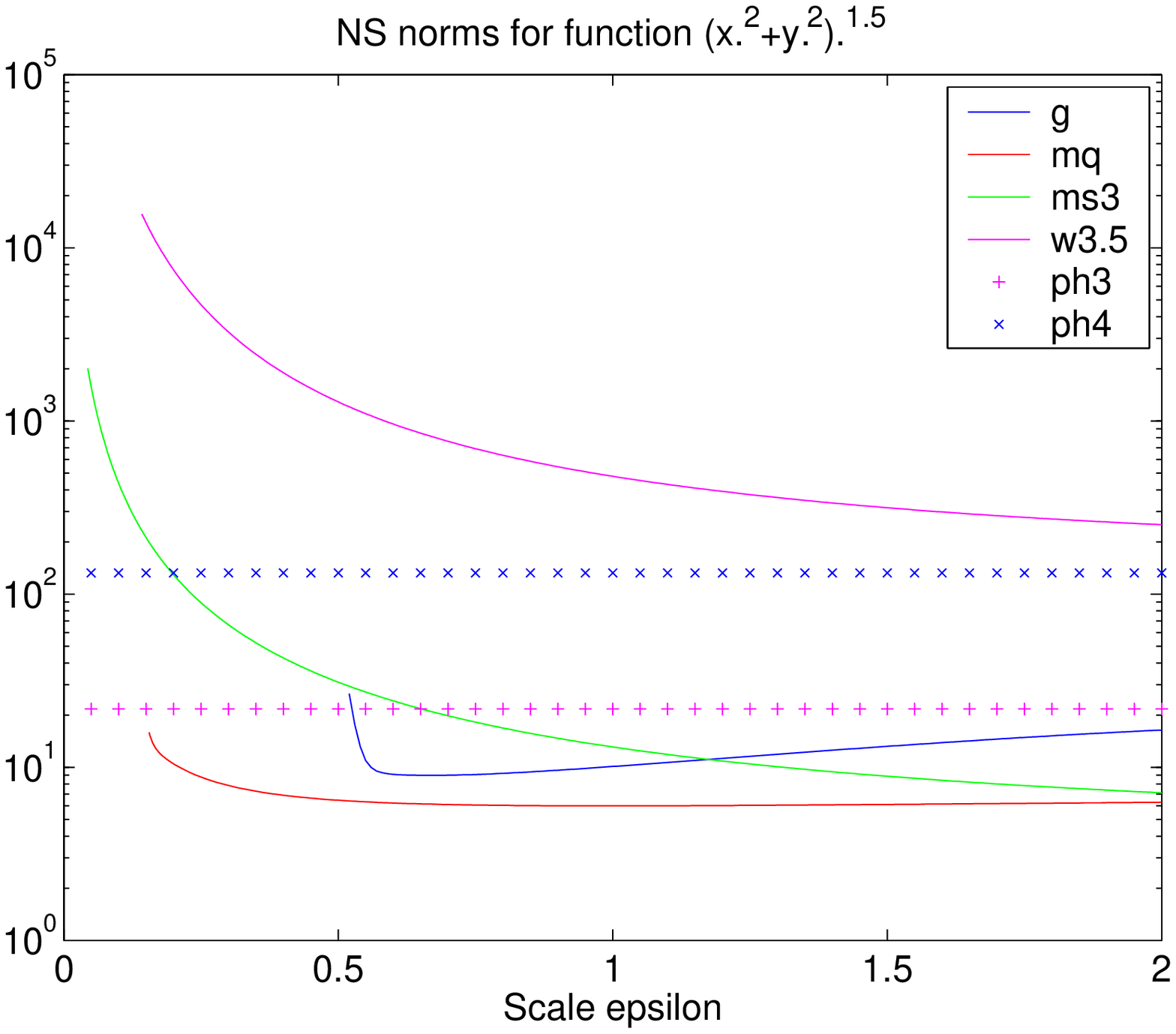}\\
\includegraphics[height=5.5cm, width=5.5cm]{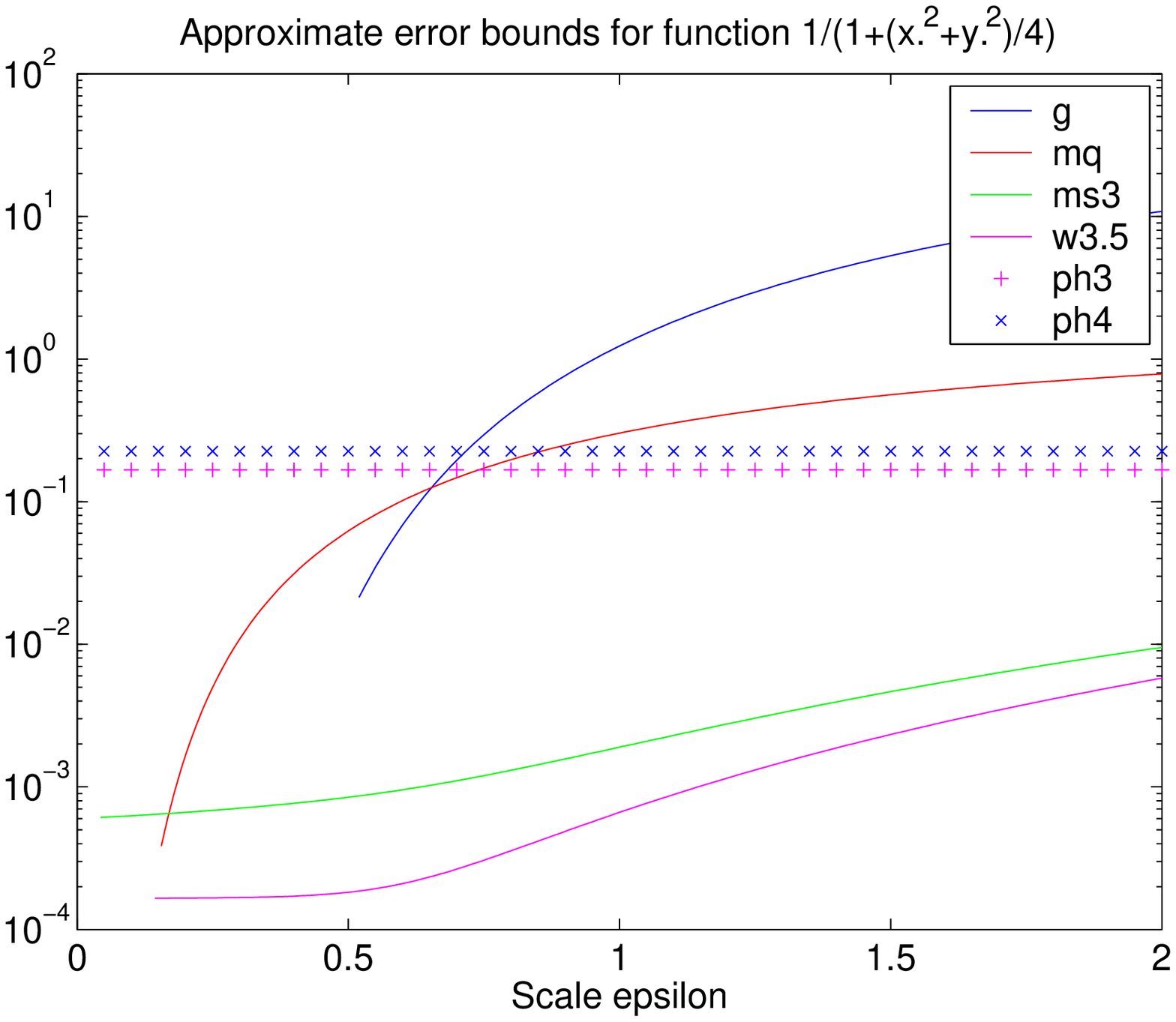}
\includegraphics[height=5.5cm, width=5.5cm]{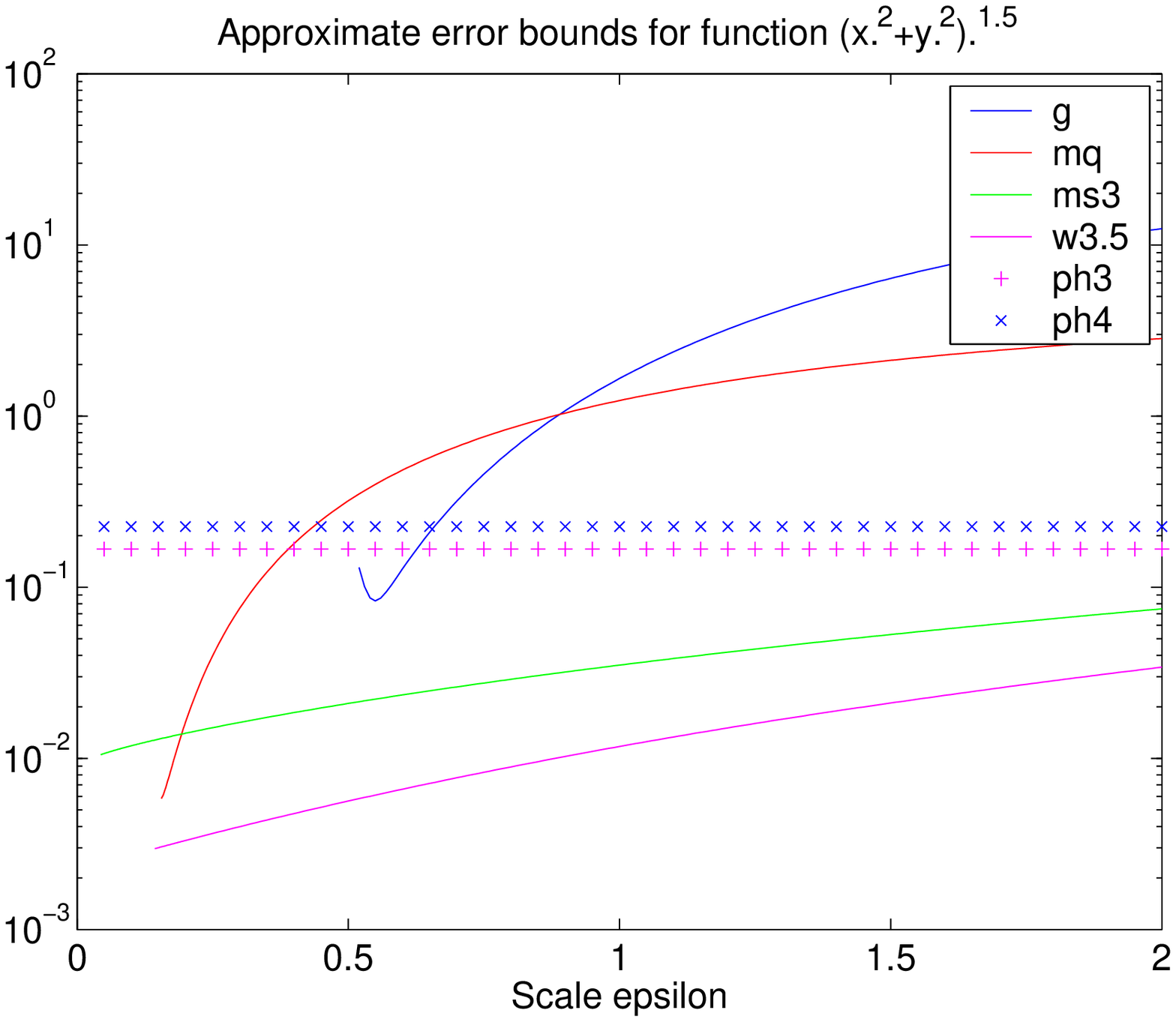}
\end{center}
\caption{Errors, native space norms, and approximate error bounds
  for two selected functions.\RSlabel{Fig28allinone}}
\end{figure}
\begin{figure}
\begin{center}
\includegraphics[height=5.5cm, width=5.5cm]{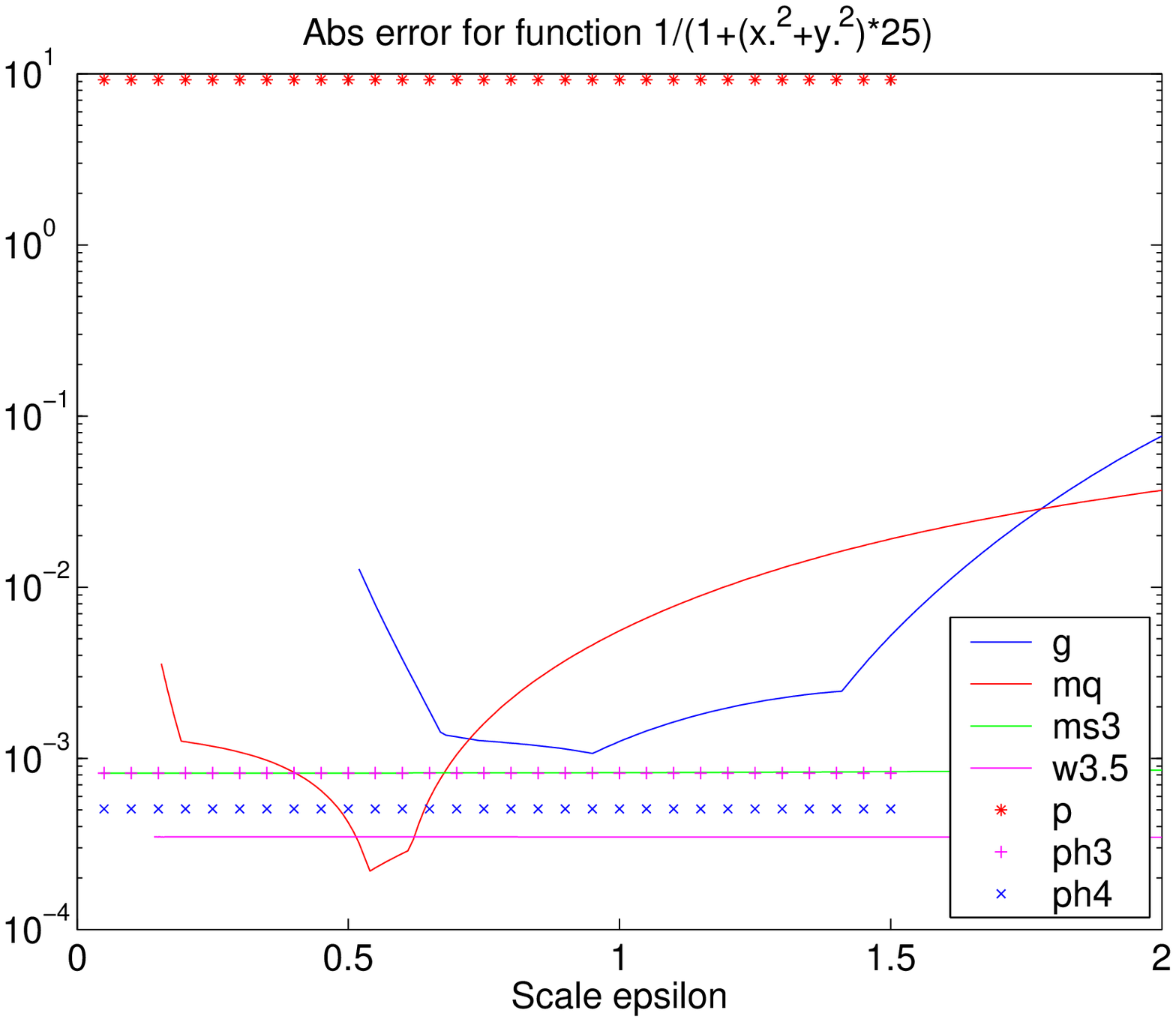}
\includegraphics[height=5.5cm, width=5.5cm]{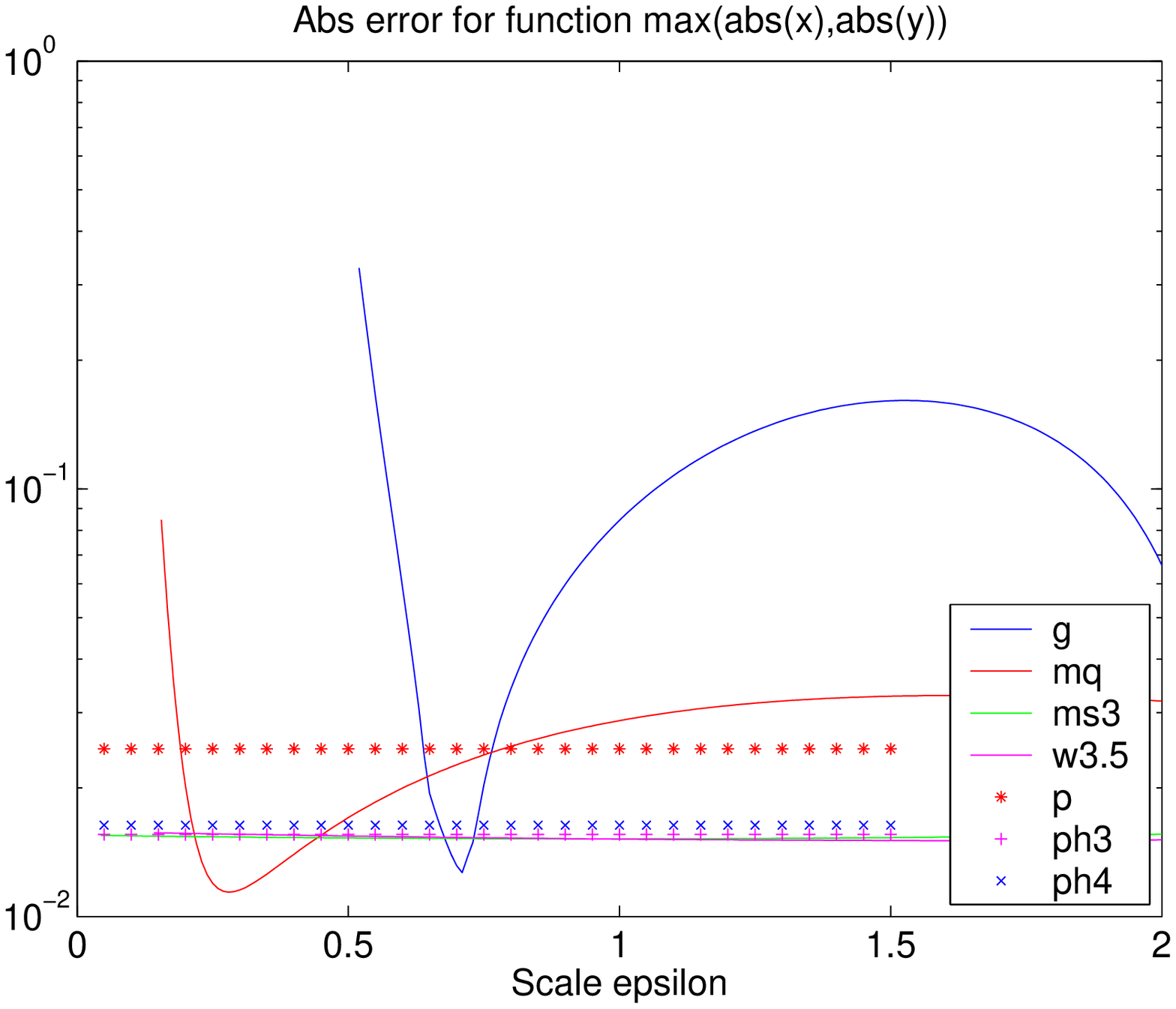}\\
\includegraphics[height=5.5cm, width=5.5cm]{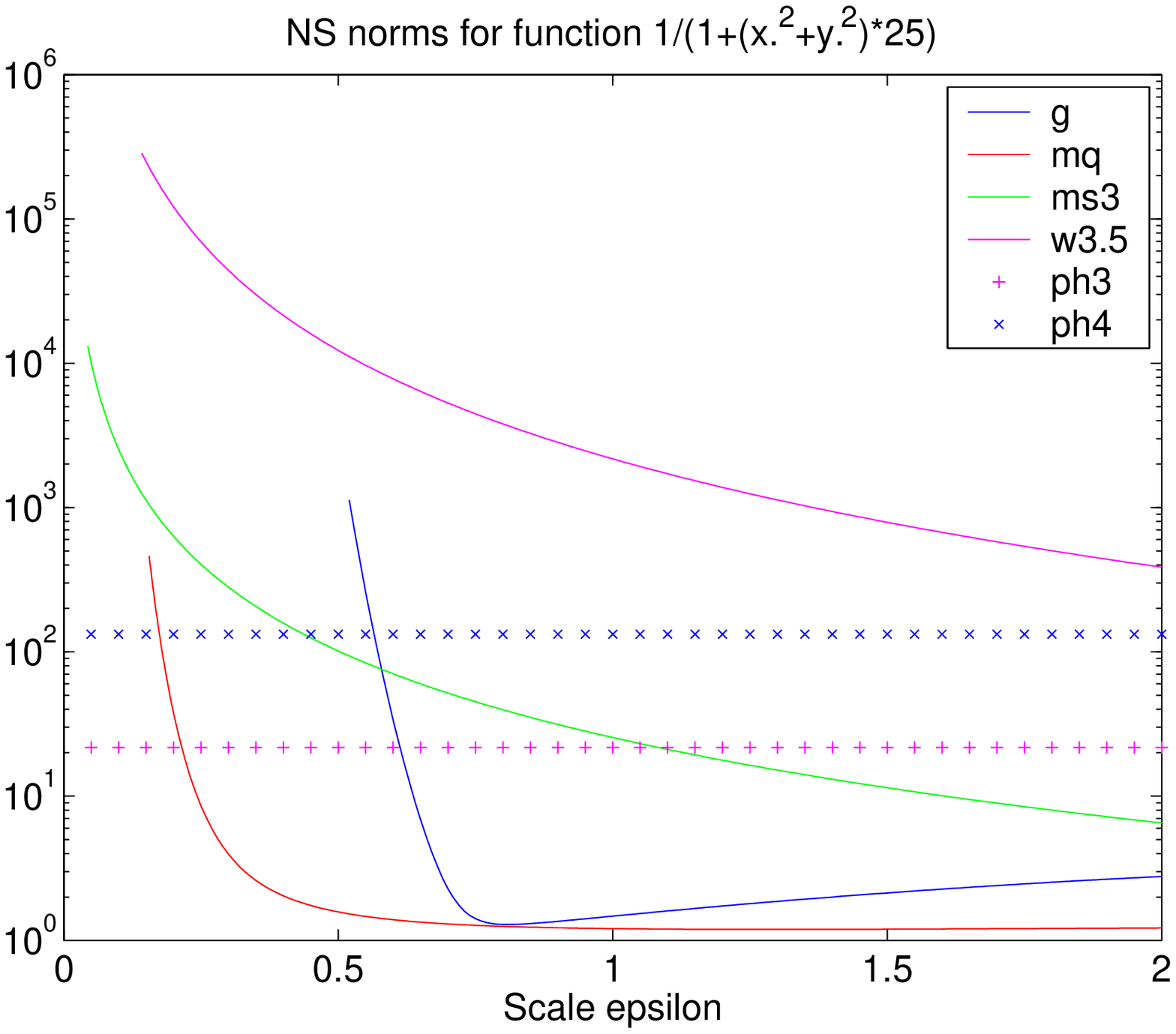}
\includegraphics[height=5.5cm, width=5.5cm]{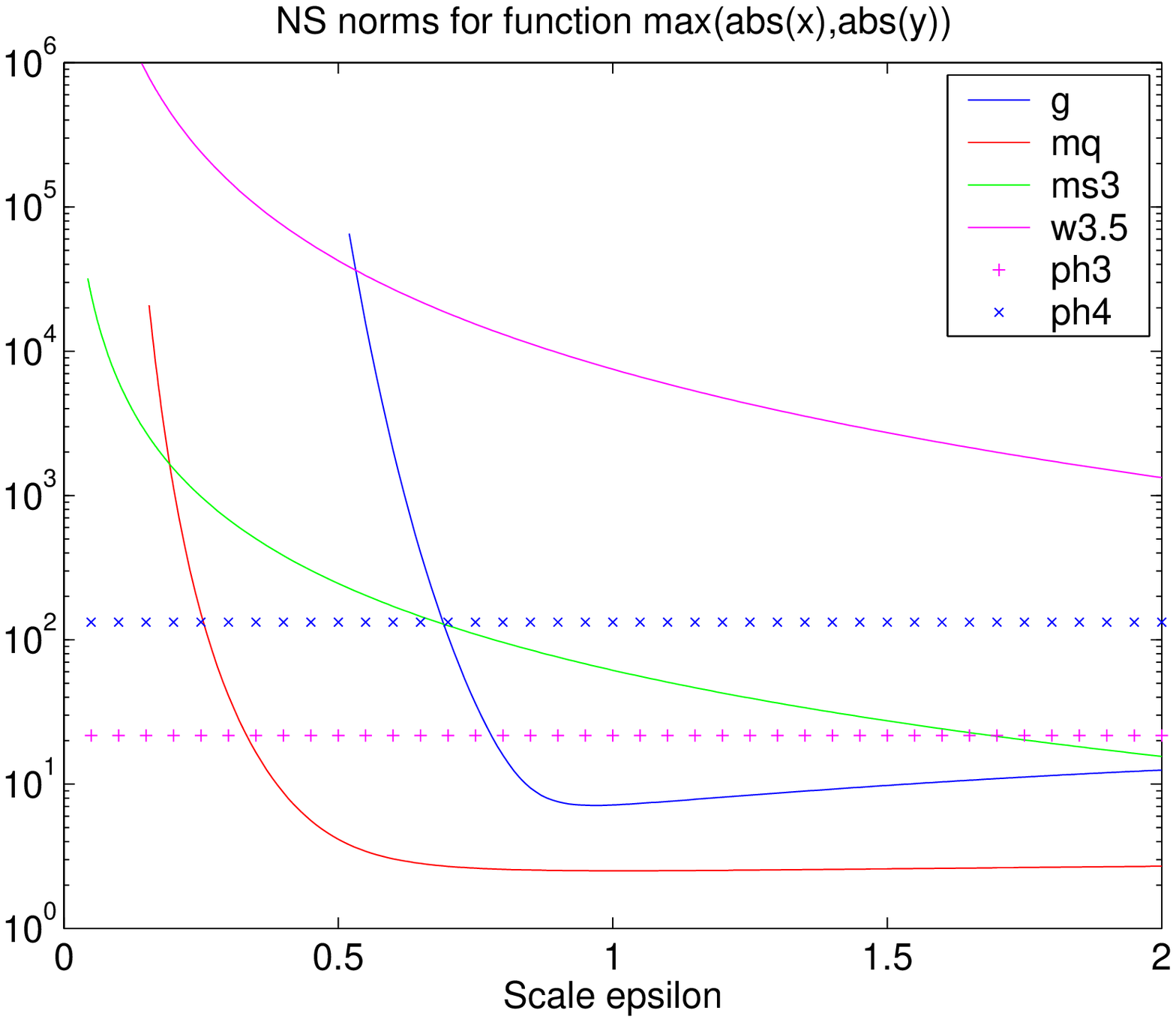}\\
\includegraphics[height=5.5cm, width=5.5cm]{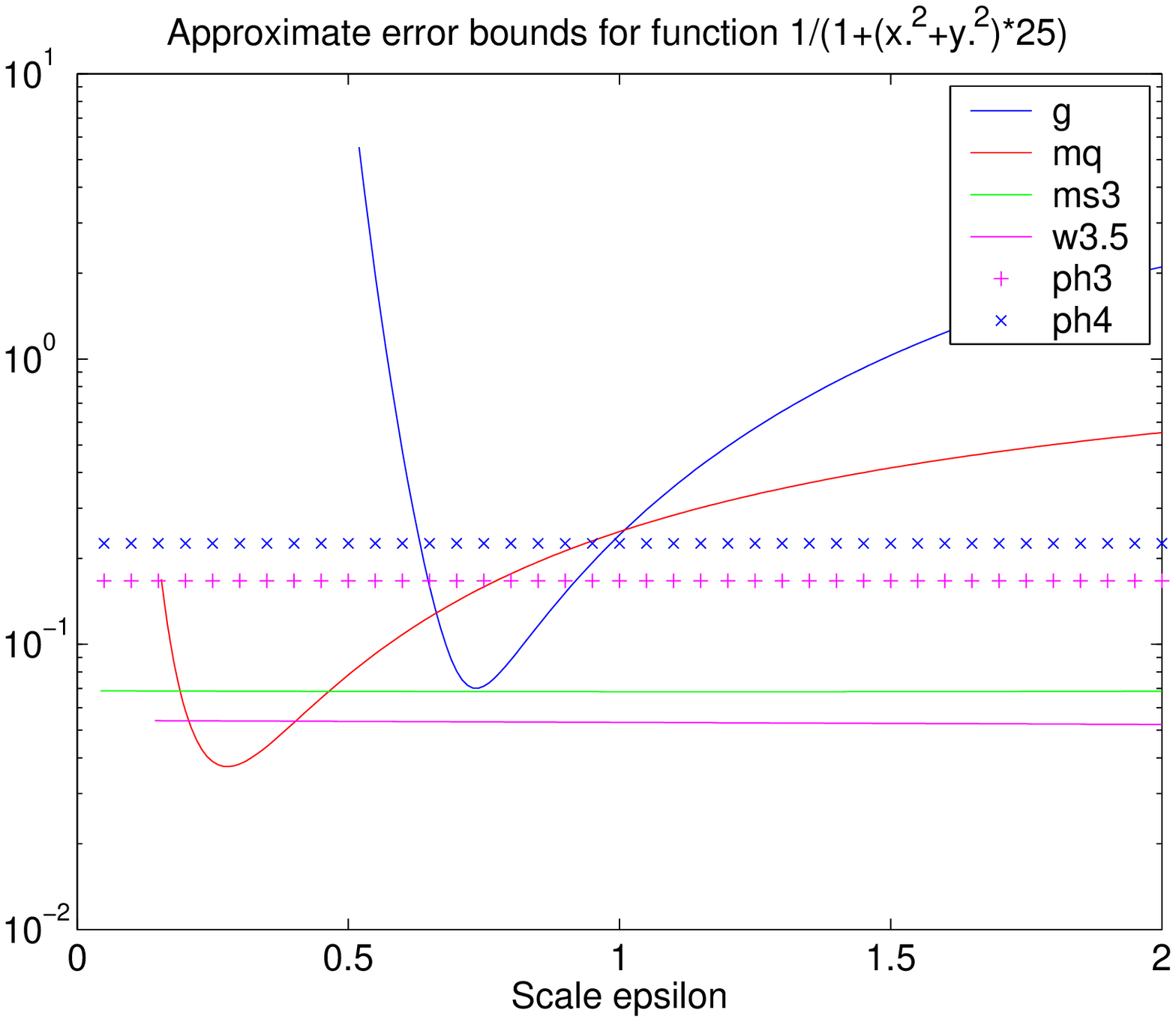}
\includegraphics[height=5.5cm, width=5.5cm]{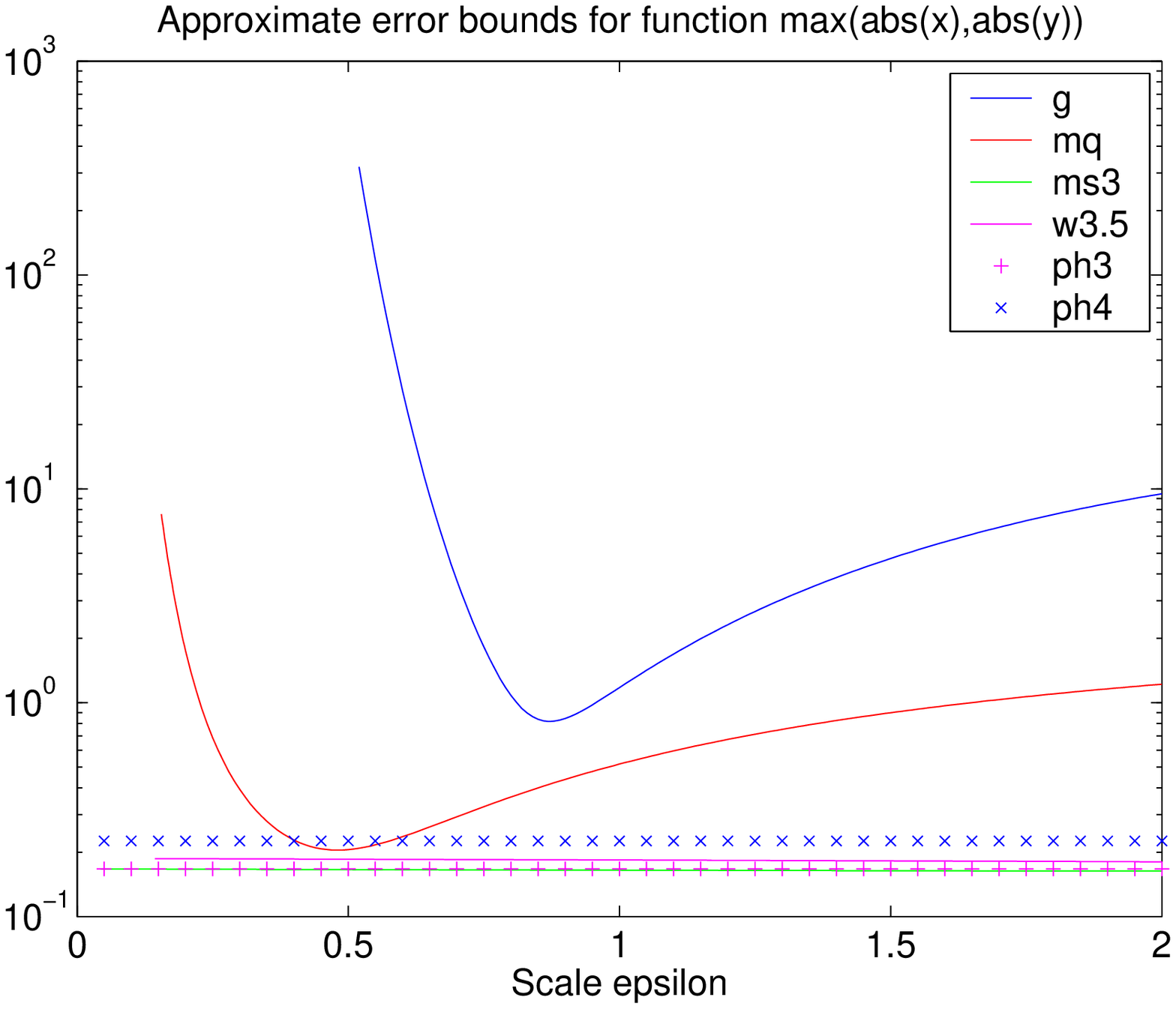}
\end{center}
\caption{Errors, native space norms, and approximate error bounds
  for two selected functions\RSlabel{Fig76allinone}}
\end{figure}
\biglf
In the flat limit, the factors in the error bound \RSref{eqfsPf}
tend towards $0\cdot \infty$ with unknown behavior. The flat limit marks
in Figure \RSref{FigscAEB441} are the error bounds for the polyharmonic case,
and not attained as the limits of the curves, unlike Figure
\RSref{FigscErr441}. The conclusion is that the error bounds are useless when
traced to the flat limit. 
%****************************************************************
\section{Error Expansions}\RSlabel{SecEE}
We now look at cases with infinite smoothness. There are cases where
a positive scale $\epsilon$ gives smaller errors than the flat
limit, and we want to see under which conditions this works.
To this end, we have to expand the interpolation error into a power series
in $\epsilon$ and see how it behaves.
\biglf
Consider the interpolation problem on a set $X$ of
$n$ points $x_1,\ldots,x_{n}$.
as in Section \ref{SecIntro}. 
For the Gaussian kernel \cite{schaback:2005-2} and other analytic radial basis
functions \cite{larsson-fornberg:2005-1},
the flat limit for $\epsilon\to 0$ exists and is a polynomial,
if certain nondegeneracy
conditions hold \cite[I-III, p. 110]{larsson-fornberg:2005-1}.
In addition, then
% by  \cite{larsson-fornberg:2005-1}
there is an expansion
\bql{eqsexperr}
s_{f,X,\Phi_\epsilon}(x)=p_{0,f}(x)+\ep{2}p_{2,f}(x)+\ep{4}p_{4,f}(x)+\cdots,
\eq
of the interpolant, where $p_{0,f}$ is the polynomial interpolant in the flat limit,
and where
the other $p_{2j,f}$ are polynomials vanishing on all data sites.
They are linear in the data of $f$ on $X$, i.e. they have the form
\bql{eqpexp}
p_{2j,f}=\displaystyle{\sum_{i=1}^{n}f(x_i)p_{2j,i}}  
\eq
with polynomials $p_{2j,i}$ that are only dependent on $X$ and the kernel
$\Phi$. If $p_{0,f}$ has degree $K$, then $p_{2j,f}$ has degree $K+2j$
\cite{larsson-fornberg:2005-1}.
%\blue{Citation?}
%The expansion can also
%be obtained by corresponding expansions of the Lagrangians into
%$$
%u_{n,\epsilon}(x)=\displaystyle{\sum_{j=0}^\infty\epsilon^{2j}p_{2j,n}(x),\;1\l%eq
%n\leq N} 
%$$
%by techniques of \RScite{schaback:2005-2,schaback:2008-2}
%and the nondegeneracy conditions
%imposed there.
%The notation $p_{2j,n}$ does not imply that the polynomials are
%even or have degree $2j$. They are just linked to the power $\epsilon^{2j}$ and
%the point $x_n$.
\biglf
When going towards scales that yield a smaller error than the flat limit
for given $f$ and $X$, consider
$$
\sigma_f(x)=\hbox{ sgn }(p_{0,f}(x)-f(x))
$$
and write the error as a perturbation of the flat limit error as 
$$
\sigma_f(x)\left( s_{f,X,\Phi_\epsilon}(x)-f(x)\right)
= |p_{0,f}(x)-f(x)|+\displaystyle{\sigma_f(x)\sum_{j=k}^\infty \epsilon^{2j}
\sum_{i=1}^{{n}}f(x_{{i}})p_{2j,{{i}}}(x)}. 
$$
Some $\epsilon>0$ can only be an improvement over the flat limit, if
\bql{eqsefp}
\displaystyle{\sigma_f(x)\sum_{j=k}^\infty \epsilon^{2j}
\sum_{{i}=1}^{{n}}f(x_{{i}})p_{2j,{{i}}}(x)}<0
\eq
on all extremal points $x$ of $f-p_{0,f}$. This might be true even for values of
$\epsilon$ that are not extremely small.  If we focus on very small $\epsilon$
and take the smallest $k$ such that $p_{2k,f}$ is nonzero, 
we need
$$
\displaystyle{\sigma_f(x)
\sum_{{{i}}=1}^{{n}}f(x_{{i}})p_{2k,{{i}}}(x)}<0.
$$
Both criteria are rather hard to check and analyze.
\biglf
Users may think that the data points $X$  and the function values
$y_j=f(x_j)$ already determine whether the flat limit is optimal or not,
and that this suffices to estimate a good scale. But this is not true.
Even if the interpolation data are fixed,
there are functions that can lead to any possible outcome.
\begin{Theorem}\RSlabel{TheAllpossible}
  For all choices of point sets $X=\{x_1,\ldots,x_{{n}}\}$ and function values
  $y_1,\ldots,y_{n}$ it cannot be decided whether the flat limit is
  error-optimal or not.
  More precisely, there are functions $f$ that attain the values
  $y_j=f(x_j),\;1\leq j\leq {n}$ and have an optimal flat limit or not.
\end{Theorem} 
\begin{proof}
Let the point set $X$ and the kernel $\Phi$ be fixed. Furthermore, we
consider fixed data values $y_j$ at the $x_j\in X$ and let the functions $f$
with $f(x_j)=y_j,\;x_j\in X$ vary. For each such function, the polynomials
$p_{2k,f}$ are the same, and we denote them now by $p_{2k,X,Y}$
to separate the set  $X$ of data locations
  from the set $Y$ of function values. Now let $g$ be
any function that vanishes on $X$ and
consider the function $f_g(x)=p_{0,X,Y}(x)-g(x)$.
Then
$$
\begin{array}{rcl}
  s_{f_g,X,\Phi_\epsilon}(x)-f_g(x)
  &=& s_{p_{0,X,Y},X,\Phi_\epsilon}(x)-p_{0,X,Y}(x)+g(x)\\
  &=& g(x)+ \displaystyle{\sum_{j=k}^\infty
    \epsilon^{2j}p_{2j,X,Y}(x)}\\
\end{array}
$$
shows that all functions $g$ that vanish on $X$ can arise as error functions
in the flat limit. Let $k(X,Y,g)$ be the smallest $k\geq 1$
such that $p_{2k,X,Y}$ does not
vanish on the extremal points of $g$. Then
$$
\hbox{sgn }p_{2k(X,Y,g),X,Y}(x)=- \hbox{sgn }g(x) 
$$
should hold on the extremal points of $g$ if the flat limit is 
not optimal for small $\epsilon$.
\biglf
This is the situation when $X,Y$, and $g$ are prescribed.  But one may
take $k(X,Y)$ to be the smallest $k\geq 1$ such that $p_{2k,X,Y}$ does not
vanish, and then fix $g=-p_{2k,X,Y}$ to arrive at the same situation, while the
choice $g=+p_{2k,X,Y}$ leads to local optimality of the flat limit. 
\end{proof}
The consequence is that users must spare data points for error evaluation,
if they want to get anywhere with a numerical estimation of scales.
This is part of
the Leave-one-out-cross-validation (LOOCV) technique,
see \RScite{fasshauer-zhang:2007-1}.
%%%%%%%%%%%%%%%%%%%%%%%%%%%%
%****************************************************************
\section{Conclusions and Open Problems}\RSlabel{SecCaOP}
Choosing ``good'' scales for kernel-based interpolation remains a
challenge, both in theory and practice. A few things could be
contributed here, but there are many open issues. 
\biglf
Sobolev spaces on $\R^d$ form a scale of spaces, and each function
belonging to one of these spaces belongs to all of them.
When minimizing the norm of a given function over all scales
of spaces, there is a unique minimum, defining a ``natural'' scale
of the function. But this scale is hard to estimate, and it is
only of limited practical  value for finding
error-optimal scales, because it minimizes a function norm,
not an interpolation error norm. However, the ``natural scale''
may be useful elsewhere. 
\biglf
If the error of an interpolant cannot be evaluated because there are
no additional data, one might look at the minimum
of the standard error bound, because it does not need additional data.
But, in contrast to the natural scale, there is no theoretical guarantee for
a nonzero scale that minimizes the error. In particular, the flat limit may
be the zero scale that yields the optimal error. There is some
experimental support for the strategy to calculate the flat limit case
in parallel to the standard scale-dependent solution. In a variety of cases,
the flat limit error comes out small enough to let the
unstable hunt for an error-optimal scale be not worth while. Criteria
for this case would be helpful. 
\biglf
In the case of analytic kernels like the Gaussian or the multiquadrics,
changes of scales imply  nontrivial changes of the underlying
Reproducing Kernel Hilbert space, making scale changes a hazardous
issue even in theory. But this is no restriction if one focuses
on weak error norms. In particular, one can expand the pointwise
error into a series wrt. the scale, but the result does not easily yield
useful criteria for ``good'' scales. These criteria should be
improved. However, the expansions
can be used for showing that for fixed interpolation data, the
unknown function on the other points can be chosen to
yield any possible error function. In addition, the flat limit may be optimal
or not, just by changing the function outside the given interpolation data.
\biglf
Readers are encouraged to extend the above situations, and this paper hopefully
serves as a starting point. 
%****************************************************************
\bibliographystyle{plain}

%\bibliography{RSbib}  % to be replaced by the bbl file

\begin{thebibliography}{}

\end{thebibliography}


\begin{thebibliography}{10}
\bibitem{baxter-brummelhuis:2022-1}
  B. Baxter and R. Brummelhuis.
  \newblock Convergence Estimates for Stationary Radial Basis Function
    Interpolation and for Semi-discrete Collocation-Schemes.
    \newblock {\em J. Fourier Anal. Appl.} 28, 53 (2022)
    % https://doi.org/10.1007/s00041-022-09945-3}
\bibitem{buhmann:1989-2}
M.D. Buhmann.
\newblock {\em Multivariable interpolation using radial basis functions}.
\newblock PhD thesis, University of Cambridge, 1989.
\bibitem{buhmann:1998-1}
M.D. Buhmann.
\newblock Radial functions on compact support.
\newblock {\em Proceedings of the Edinburgh Mathematical Society}, 41:33--46,
  1998.

\bibitem{buhmann:2003-1}
M.D. Buhmann.
\newblock {\em Radial Basis Functions, Theory and Implementations}.
\newblock Cambridge University Press, Cambridge,UK, 2003.

\bibitem{davydov-schaback:2019-1}
O.~Davydov and R.~Schaback.
\newblock Optimal stencils in {S}obolev spaces.
\newblock {\em IMA Journal of Numerical Analysis}, 39:398--422, 2019.

\bibitem{demarchi-schaback:2010-1}
St. De~Marchi and R.~Schaback.
\newblock Stability of kernel-based interpolation.
\newblock {\em Adv. in Comp. Math.}, 32:155--161, 2010.

\bibitem{driscoll-fornberg:2002-1}
T.A. Driscoll and B.~Fornberg.
\newblock Interpolation in the limit of increasingly flat radial basis
  functions.
\newblock {\em Comput. Math. Appl.}, 43:413--422, 2002.

\bibitem{fasshauer-mccourt:2015-1}
G.~Fasshauer and M.~McCourt.
\newblock {\em Kernel-based Approximation Methods using MATLAB}, volume~19 of
  {\em Interdisciplinary Mathematical Sciences}.
\newblock World Scientific, Singapore, 2015.

\bibitem{fasshauer-zhang:2007-1}
G.~Fasshauer and J.G. Zhang.
\newblock On choosing optimal shape parameters for {RBF} approximation.
\newblock {\em Numerical Algorithms}, 45:345--368, 2007.

\bibitem{fornberg-et-al:2011-1}
B.~Fornberg, E.~Larsson, and N.~Flyer.
\newblock Stable computations with {G}aussian radial basis functions.
\newblock {\em SIAM J. Sci. Comput.}, 33(2):869--892, 2011.

\bibitem{fornberg-et-al:2013-1}
B.~Fornberg, E.~Lehto, and C.~Powell.
\newblock Stable calculation of {G}aussian-based {RBF-FD} stencils.
\newblock {\em Computers and Mathematics with Applications}, 65:627--637, 2013.

\bibitem{fornberg-wright:2004-1}
B.~Fornberg and G.~Wright.
\newblock Stable computation of multiquadric interpolants for all values of the
  shape parameter.
\newblock {\em Comput. Math. Appl.}, 48(5-6):853--867, 2004.

\bibitem{fornberg-et-al:2004-1}
B.~Fornberg, G.~Wright, and E.~Larsson.
\newblock Some observations regarding interpolants in the limit of flat radial
  basis functions.
\newblock {\em Computers \& Mathematics with Applications}, 47:37--55, 2004.
\newblock doi:10.1016/S0898-1221(04)90004-1.

\bibitem{kansa-carlson:1992-1}
E.~J. Kansa and R.E. Carlson.
\newblock Improved accuracy of multiquadric interpolation using variable shape
  parameters.
\newblock {\em Comput.\ Math.\ Appl.}, 24:99--120, 1992.

\bibitem{larsson-fornberg:2005-1}
Elisabeth Larsson and Bengt Fornberg.
\newblock Theoretical and computational aspects of multivariate interpolation
  with increasingly flat radial basis functions.
\newblock {\em Comput.\ Math.\ Appl.}, 49:103--130, 2005.

\bibitem{lee-et-al:2005-1}
Y.J. Lee, G.J. Yoon, and J.~Yoon.
\newblock Convergence property of increasingly flat radial basis function
  interpolation to polynomial interpolation.
\newblock preprint, 2005.

\bibitem{schaback:1995-1}
R.~Schaback.
\newblock Error estimates and condition numbers for radial basis function
  interpolation.
\newblock {\em Advances in Computational Mathematics}, 3:251--264, 1995.

\bibitem{schaback:2005-2}
R.~Schaback.
\newblock Multivariate interpolation by polynomials and radial basis functions.
\newblock {\em Constructive Approximation}, 21:293--317, 2005.

\bibitem{schaback:2008-2}
R.~Schaback.
\newblock Limit problems for interpolation by analytic radial basis functions.
\newblock {\em J. Comp. Appl. Math.}, 212:127--149, 2008.

\bibitem{scheuerer-et-al:2013-1}
M.~Scheuerer, M.~Schlather, and R.~Schaback.
\newblock Interpolation of spatial data - a stochastic or a deterministic
  problem?
\newblock {\em European Journal of Applied Mathematics}, 24:601--629, 2013.

\bibitem{song-et-al:2012}
G.~Song, J.~Riddle, G.E. Fasshauer, and F.J. Hickernell.
\newblock Multivariate interpolation with increasingly flat radial basis
  functions of finite smoothness.
\newblock {\em Adv. Comp. Math.}, 36:485--501, 2012.

\bibitem{wendland:1995-1}
H.~Wendland.
\newblock Piecewise polynomial, positive definite and compactly supported
  radial functions of minimal degree.
\newblock {\em Advances in Computational Mathematics}, 4:389--396, 1995.

\bibitem{wendland:2005-1}
H.~Wendland.
\newblock {\em Scattered Data Approximation}.
\newblock Cambridge University Press, Cambridge,UK, 2005.

\end{thebibliography}
\end{document}